\documentclass[a4paper,12pt,numbers,sort&compress]{article}

\pdfoutput=1

\usepackage[paper=letterpaper,margin=1in]{geometry}

\usepackage{amsmath,amssymb,amsfonts,epsfig,cite,setspace,bigstrut,longtable,array,url,color}

\usepackage{graphicx}
 \usepackage[textwidth=3.4cm,textsize=10pt]{todonotes}

\usepackage[singlelinecheck=false,justification=justified]{caption}

\usepackage{hyperref}
\usepackage{pdfpages,lipsum}

\renewcommand{\baselinestretch}{1.028}
\newcommand{\refpart}[1]{{\it (#1)}}  

\begin{document}


\vskip 0.25in

\newcommand{\sref}[1]{\S~\ref{#1}}
\newcommand{\nn}{\nonumber}
\newcommand{\tr}{\mathop{\rm Tr}}
\newcommand{\comment}[1]{}
\newcommand{\cW}{{\cal W}}
\newcommand{\cN}{{\cal N}}
\newcommand{\cQ}{{\cal Q}}
\newcommand{\cH}{{\cal H}}
\newcommand{\cK}{{\cal K}}
\newcommand{\cZ}{{\cal Z}}
\newcommand{\cO}{{\cal O}}
\newcommand{\cA}{{\cal A}}
\newcommand{\cB}{{\cal B}}
\newcommand{\cC}{{\cal C}}
\newcommand{\cD}{{\cal D}}
\newcommand{\cE}{{\cal E}}
\newcommand{\cF}{{\cal F}}
\newcommand{\cG}{{\cal G}}
\newcommand{\cM}{{\cal M}}
\newcommand{\cX}{{\cal X}}
\newcommand{\IA}{\mathbb{A}}
\newcommand{\IP}{\mathbb{P}}
\newcommand{\IQ}{\mathbb{Q}}
\newcommand{\IH}{\mathbb{H}}
\newcommand{\IR}{\mathbb{R}}
\newcommand{\IC}{\mathbb{C}}
\newcommand{\IF}{\mathbb{F}}
\newcommand{\IV}{\mathbb{V}}
\newcommand{\II}{\mathbb{I}}
\newcommand{\IZ}{\mathbb{Z}}
\newcommand{\re}{{\rm~Re}}
\newcommand{\im}{{\rm~Im}}

\newcommand{\ol}{\overline}
\newcommand{\diff}{\partial}
\newcommand{\dbar}{\ol{\partial}}

\newcommand{\tmat}[1]{{\tiny \left(\begin{matrix} #1 \end{matrix}\right)}}
\newcommand{\mat}[1]{\left(\begin{matrix} #1 \end{matrix}\right)}

\newcommand{\hpg}[5]{{}_{#1}\mbox{\rm F}_{\!#2}\!
  \left(\left.{#3 \atop #4}\right| #5 \right) }
\newcommand{\hpgo}[2]{{}_{#1}\mbox{\rm F}_{\!#2}}

\let\oldthebibliography=\thebibliography
\let\endoldthebibliography=\endthebibliography
\renewenvironment{thebibliography}[1]{%
\begin{oldthebibliography}{#1}%
\setlength{\parskip}{0ex}%
\setlength{\itemsep}{0ex}%
}%
{%
\end{oldthebibliography}%
}

\newtheorem{theorem}{Theorem}[section]
\newtheorem{remark}[theorem]{Remark}

\def\theequation{\thesection.\arabic{equation}}
\newcommand{\setall}{\setcounter{equation}{0}
        \setcounter{theorem}{0}}
\newcommand{\setequation}{\setcounter{equation}{0}}
\renewcommand{\thefootnote}{\fnsymbol{footnote}}

~\\
\vskip 1cm

\begin{center}
{\Large \bf Composite Genus One Belyi Maps} 
\end{center}
\medskip

\vspace{1cm}
\centerline{
{\large Raimundas Vidunas}$^1$ \ \&
{\large Yang-Hui He}$^2$
}
\vspace*{4.0ex}

\renewcommand{\baselinestretch}{0.5}
\begin{center}
{\it
  {\small \qquad \quad
    \begin{tabular}{cl}
      ${}^{1}$ 
      & Graduate School of Information Science and Technology,\\
     & Osaka University, Osaka, Japan.\\
      ${}^{2}$ 
      & Department of Mathematics, City, University of London, EC1V 0HB, UK \& \\
      & School of Physics, NanKai University, Tianjin, 300071, P.R.~China \& \\
      & Merton College, University of Oxford, OX14JD, UK\\
      &\\
      & \qquad
      {\rm rvidunas@gmail.com, hey@maths.ox.ac.uk}
    \end{tabular}
  }
}
\end{center}

\renewcommand{\baselinestretch}{1.028}

\vspace*{7.0ex}
\centerline{\textbf{Abstract}} \bigskip

Motivated by a demand for explicit genus 1 Belyi maps from theoretical physics, 
we give an efficient method of explicitly computing genus one Belyi maps 
by (1) composing covering maps from elliptic curves to the 
Riemann sphere  with simpler (univariate) genus zero Belyi maps,
as well as by (2) composing further with isogenies between elliptic curves.
This gives many new explicit dessins on the doubly periodic plane, including several which have been realized in the physics literature as so-called brane-tilings in the context of quiver gauge theories.

\thispagestyle{empty}

\newpage

\tableofcontents

\section{Introduction and Summary}
\setall

One of the most profound pieces of modern mathematics is the so-called ``sketches'', viz., {\it l'Esquisse d'une programme}, by Grothendieck, where he visualized the theorem of Belyi \cite{belyi} as bipartite graphs on Riemann surfaces which he called {\it dessins d'enfants} \cite{leila}, or ``children's drawings'', and thence proceeded to obtain a unified framework of number theory, algebraic geometry and combinatorics.
Belyi's theorem itself is remarkably simple: a compact, smooth Riemann surface $\Sigma$ admits an algebraic model over $\overline{\IQ}$ IFF there is a rational sujective map $\beta : \Sigma \to \IP^1$ ramified at exactly 3 points (which can be chosen to be $\{0,1,\infty\}$.
Take a continuous path $\gamma \subset \IP^1$ with endpoints at $0$ and $1$, say.
Then the pre-image $\beta^{-1}(\gamma)$ is a graph embedded in $\Sigma$ which is bi-partite by having all pre-images of 0 coloured, say, black, and that of 1, white;
this is the dessin.
Finding the explicit map $\beta$ and algebraic equation for $\Sigma$ (as, for instance, a hyper-elliptic curve), on the other hand, is often a formidable task and is itself a long programme \cite{catalog,hv,jones,kmsv,sv,vk,vf,vhpg,vhpgell,vHVHeun}.

In parallel, bipartite graphs on Riemann surfaces, especially of genus 1 -- called ``dimer models'' or {\it brane tilings} -- have been a long subject in physics within the context of contructing supersymmetric gauge theories from branes in string theory \cite{Franco:2005rj,Franco:2005sm} (cf.~\cite{He:2012js} for a rapid review), which has been for a number of years having a fruitful dialogue with the algebraic geometry (especially Calabi-Yau varieties) and quiver representation community (cf.~a rather comprehensive treatment in \cite{Broom}).
The inter-play between the physics of the gauge theory, the algebraic geometry of the underlying toric Calabi-Yau geometry probed by the branes and the representation theory of the quiver structure which encodes both has been inspiring many new directions of research 

It is therefore natural to consider these quiver gauge theories in the guise of brane tilings as dessins on the elliptic curve;
this was undertaken recently and a subsequent programme in studying the gauge theory-dessin-Calabi-Yau correspondence
was also established \cite{Jejjala:2010vb,Hanany:2011ra,Hanany:2011bs},
and the significance of proceeding beyond genus 1, also pointed out \cite{Bose:2014lea,Cremonesi:2013aba}.
In brief, let $D$ be a dessin on a genus $g$ Riemann surface $\Sigma_g$, the graph dual to this is a quiver with relations coming from the Jacobian of a polynomial $W$; in turn, this is a supersymmetric gauge theory in $3+1$-dimensions whose vacuum moduli space is an affine toric Calabi-Yau variety of dimension $2g+1$.
Hence, the case of $g=1$ gives the particularly interesting case of Calabi-Yau threefolds which are central to string theory.

In light of the foregoing discussions, systematic investigations of genus 1 Belyi maps, i.e., dessins d'enfants on an elliptic curve, is thus an important endeavour.
Some of these relevant genus 1 dessins, presented as brane tilings on the doubly periodic plane, are recently classified in \cite{Davey:2009bp}.
However, for only very few of these theories have the corresponding Belyi pair 
of $(\Sigma, \beta)$ been explicitly constructed.
This is partly because the computation of genus 1 maps is significantly harder than in the genus 0 case.

One easy way to construct genus 1 Belyi maps is to look for compositions of $\varphi_0$
with a quadratic covering $\pi:E\to \IP^1$ that branches over 4 points.
If those 4 points are in the 3 branching fibers of $\varphi_0$,
the composite map $\pi\circ\varphi_0$ will be a Belyi map of genus 1.
This simple observation will be our starting point, which, together with isogenies on elliptic curves, will turn out to be a powerful tool to produce a host of new genus 1 dessins.
These two basic constructions were suggested 
in \cite{Hanany:2011ra}, 
as Belyi maps ``dependent on $x$ only" and orbifolds \cite[\S 5]{Jejjala:2010vb}, respectively.

The plan of the paper is as follows.
In \S\ref{sec:comp} we find compositions of this kind for examples within the catalogue \cite{Davey:2009bp} and realize them, this gives many new dessin realization of the physical model of brane-tilings.
Mathematically, we arrive at a multitude of dessins on the torus of degrees 10, 12, up to 18. 
Next, in \S \ref{sec:isogenies}, we also consider construction of genus 1 Belyi maps by composing known genus 1 Belyi maps with isogenies between elliptic curves. In particular, we construct isogenies to the curve $Y^2=X^3-X$ and the composite Belyi maps with the square lattice dessins in by reverting to reduction to genus 0 dessins.
Finally, in \S \ref{sec:othercomp}, we discuss composition of genus 0 Belyi maps with higher degree coverings $E\to\IP^1$, notably using a paramtric degree 3 covering and its specialization to the highly symmetric $j=0$ elliptic curve $y^3 = z(z-1)$.
We conclude with a the physical motivations for our theory in \S\ref{sec:conc}, 
 emphasizing the origin of many of our dessins as supersymmetric gauge theories from string theory, as well as prospects of future work and upcoming papers.
In the appendix, we draw the genus 1 counterparts of the genus 0 dessins presented in the main text, after the various appropriate compositions; this serves as a partial catalogue of the dessin realization of the bipartite tilings in our context.


\subsection*{Nomenclature}
\begin{itemize}
\item We use the multiplicative notation
  $[a_1^{p_1}\cdots a_{k_0}^{p_{k_0}} / b_1^{q_1}\cdots b_{k_1}^{q_{k_1}} /
  c_\infty^{r_1}\cdots c_{k_0}^{r_{k_\infty}}]$ for the branching passport
  of Belyi maps. This gives the branching 
  indices, with repetition written as exponents, 
  of the pre-images of $\infty,0,1$ respectively. 
\item We consider genus 1 curves
in a Weierstra\ss\ form $Y^2 = X^3+\ldots$ as {\em elliptic curves} \cite{EMod}
in the standard way, by considering the infinite point as the origin of the group law.
In \S \ref{sec:isogenies} and \S \ref{sec:isogj0}, 
genus 1 curves in alternative forms are considered as elliptic curves by provided isomorphisms.
\item The Klein J-invariant is denoted as $j$ and we adhere to the normalization that for $y^2 = x^3 -x$, $j=1728$  (not $j=1$).
\item The dessins d'enfant of genus 1 Belyi maps 
are most illustratively represented as doubly periodic tilings. 
In this article, the black and white nodes represent the points above $\infty,0$
(respectively), and the cells represent the points above $1$.
\item If the $j$-invariant of a genus 1 curve is real, we draw 
its fundamental region in a doubly periodic
tiling as a rectangle (if the curve can have two real components)
or as a rhombus (if it can have one real component). 
If additionally the Belyi map is defined over $\IR$, the drawn tiling will have reflexion symmetries representing the complex conjugation:  
two reflexions parallel to a side of a rectangular fundamental region, 
or a reflexion parallel to a diagonal of a rhombus. 
In this way, our pictures reflect basic symmetries of the elliptic curves and the dessins.
In particular, a square fundamental domain signifies an elliptic curve with $j=1728$.
\item We let $\omega_n=\exp(2\pi i/n)$ denote the primitive $n$-th root of unity.
\item The Gau\ss{} hypergeometric function is defined via the Euler integral: 
  \[
  \hpg{2}{1}{a,b}{c}{z} = \frac{\Gamma(c)}{\Gamma(b)\Gamma(c-b)}
  \int_0^1 t^{b-1}(1-t)^{c-b-1}(1 - tz)^{-a} dt.
  \]
It satisfies the hypergeometric differential equation 
\begin{equation} \label{eq:HGE}
z\,(1-z)\,\frac{d^2y(z)}{dz^2}+
\left(c-(a+b+1)\,z\right)\,\frac{dy(z)}{dz}-a\,b\,y(z)=0.
\end{equation}  
This is a canonical Fuchsian equation with 3 singularities, namely $z=0$, $z=1$, $z=\infty$.
Alternative symmetric parameters are the {\em local exponent differences} $1-c,c-a-b,a-b$
at the three singular points.
\end{itemize}

\section{Compositions with Genus 0 Maps}
\label{sec:comp}\setall
An easy way to get a genus 1 Belyi map is to compose a genus 0 Belyi map $\varphi_0$  with a quadratic covering $\pi$ that branches over 4 points.
The characteristic feature is that these Belyi maps can be defined as functions in the $x$-coordinate of the elliptic curve only.
In this section, we demonstrate this simple recipe to obtain genus 1 Belyi maps from genus 0 ones and use it to immediately obtain the explicit forms of dessins which are much needed in the literature.

\subsection{Computational Details}
The genus 0 maps are obtained by using the {\sf Maple} program developed by the first author
in collaboration with Mark van Hoeij in \cite{vHVHeun}. 
The program can effectively compute Belyi maps of genus 0 if the passport is nearly regular.
Given positive integers $k,\ell,m,n$,
a Belyi map in \cite{vHVHeun} is defined to be {\em $(k,\ell,m)$-minus-$n$ regular} 
if, with exactly $n$ exceptions in total, all points above $z=1$ have branching order $k$, 
all points above $z=0$ have branching order $\ell$, 
and all points above $z=\infty$ have branching order $m$. 
Typically, these Belyi maps can pull-back the hypergeometric equation (\ref{eq:HGE})
with the local exponent differences $1/k,1/\ell,1/m$ to Fuchsian equations with $n$ singularities. 
When $n$ is small, this gives additional algebraic equations 
for the undetermined coefficients of the desired Belyi maps.
See \cite[\S 5.2]{vHVHeun}, \cite{hv} for the full system of utilized relations
of undertermined coefficients.

In \cite{vHVHeun}, 366 Galois orbits of Belyi $(k,\ell,m)$-minus-$4$ Belyi maps of degree $\le 60$
were computed. Not all computed Belyi maps of genus 0 in the present article are highly regular, 
as we have $n\le 6$. But their degree is at most 12, so they can be computed within minutes.
Nevertheless, the computed examples (especially for \S \ref{sec:square}) 
revealed a few mistakes in the {\sf Maple} program.
In particular, towers of field extensions (for the definition field)
and some formations of degenerate 
solutions were not handled well.
One applied measure was to eliminate degenerate solutions more aggressively from
an originally obtained system for undetermined coefficients.

The genus 1 curves are determined by the branching points of the quadratic covering $\pi$.
These are the points where the branching order of a genus 0 map 
is doubled (rather than duplicated).
They are determined once the passport of a relevant genus 0 map 
is decided, and that Belyi map is computed.

Simplification of computed Belyi maps (by M\"obius transformations
or isomorphisms of their genus 1 curve) is a cumbersome problem,
still not satisfactorily automatized. 
In particular, bringing an equation for a genus 1 curve to 
the most canonical Weierstra\ss\ form $y^2 = 4x^3 - g_2 x - g_3$ 
is not always practical.
Often our construction from a genus 0 map 
directly gives elliptic curves in the form $y^2=x^4+ax^3+bx^2+cx+d$.
By keeping a point at infinity, we can transform this curve to the Weierstra\ss\ form 
\begin{equation} \label{eq:ecwfrom4}
Y^2=X^3+bX^2+(ac-4d)X+(a^2-4b)d+c^2
\end{equation}
by the transformation
\begin{align} \label{eq:ec4to3}
x= -\frac{2Y+aX+2c}{4X-a^2+4b}, \quad
y=  \frac{8cY-4X^3+4(ac-4d)X+8c^2+(a^2-4b)Q}{(4X-a^2+4b)^2}
\end{align}
with $Q=aY+3X^2+2bX+a c+4d$.

\begin{remark} \rm \label{rm:nonexist}
Non-existence of genus 0 Belyi maps with a given passport can be easily proven 
by contradiction to implied pull-back transformations to non-existent
Fuchsian equations \cite[\S 5]{vf}, \cite[Remark 4.1]{vHVHeun}. 
In particular, this excludes many $(k,\ell,m)$-minus-$n$ regular maps with $n\le 2$.
For example, there are no Belyi maps with the passports $[5\,3/2^4/2^4]$,
$[4\,2/2^3/2^3]$ and  $[4^2\,2/2^3\,4/2^5]$. 
\end{remark}


\subsection{Two Genus 1 Belyi Maps of Degree 12}
\label{sec:deg12a}

\begin{figure}
\begin{picture}(460,290)
\put(-2,4){\includegraphics[width=440pt]{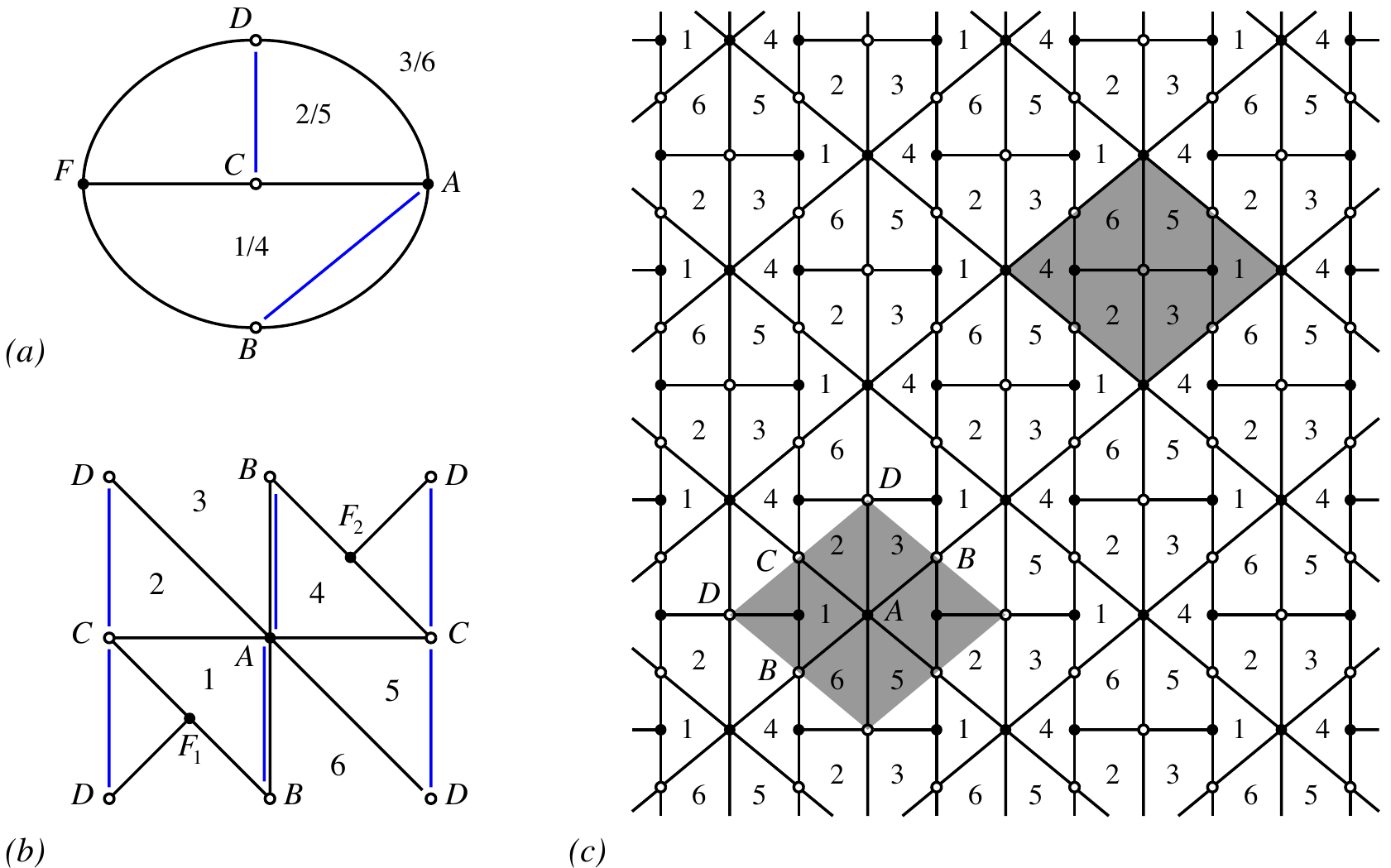}}
\end{picture}
\caption{Tiling $dP_3$(I) in 
\cite[(3.28)]{Davey:2009bp} from a genus 0 map}
\label{fig:dessin328}
\end{figure}

We start with a genus 0 Belyi map of degree 6:
\begin{equation} \label{eq:comp1}
\varphi_1 = \frac{(x^3+1)^2}{4x^3}=1+\frac{(x^3-1)^2}{4x^3}
\end{equation}
Its passport can be readily checked to be $[3^2/2^3/2^3]$ and the dessin is shown in Figure \ref{fig:dessin328}\refpart{a}.
Now, consider the $j=0$ elliptic curve $E : y^2=x^3+1$, which is a quadratic covering $E\to\IP^1$. We mark in the figure, the two cuts between points $A$, $B$, $C$, $D$ on the $\IP^1$ that are the branching points for this covering.
Most symmetrically, the cut between $A$, $B$ can be chosen to be along
the dessin segment $AB$.

Now, let $\varphi_1$ be of the same form as the above $\varphi_1$ but now as a function on $E$, and depending only on the $x$-coordinate.
We easily check that now $\varphi_1(x,y)$ a degree 12 genus 1 Belyi map with passport $[3^2\,6/4^3/2^6]$. 
The pulled-back dessin on a fundamental region of $E$ is depicted 
in Figure \ref{fig:dessin328}\refpart{b}.
The four branching points become the evenly spaced 2-torsion points on the fundamental domain.
The point $F$ and all three cells will have two preimages on $E$,
hence the double notation 1/4, 2/5, 3/6 for the cells in part \refpart{a}.

Figure \ref{fig:dessin328}\refpart{c} depicts the doubly periodic dessin,
with the gray squares representing the fundamental domain twice.
The lower domain corresponds directly to Figure \ref{fig:dessin328}\refpart{b} (after a rotation).
The upper shifted domain directly reflects the symmetries of the dessin. 
Importantly, we recognize this genus 1 dessin as \cite[(3.28)]{Davey:2009bp}.
Thus we have obtained a dessin for the so-called $dP_3$(I) theory, 
unknown to the literature thus far.


Let us continue with another example.
Consider the genus 0 Belyi map
\begin{equation} \label{eq:phi28}
\varphi_2=-\frac{(x^2-4)^3}{27x^4}=1-\frac{(x^2-1)(x^2+8)^2}{27x^4} \ ;
\end{equation}
this has the passport $[3^2/4\;2/2^21^2]$. 
Now, lift the $\IP^1$ to the elliptic curve
\begin{equation} \label{eq:ec28}
y^2=(x+1)(x-1)(x-2)
\end{equation} 
with $j$-invariant $28^3/9$.
Recasting \eqref{eq:phi28} to be a Belyi map on this elliptic curve gives a genus 1 map with the same degree 12 passport $[3^2\,6/4^3/2^6]$.

\begin{figure}[p] 
\begin{center}
\begin{picture}(382,120)
\put(-0,-8){\includegraphics[height=132pt]{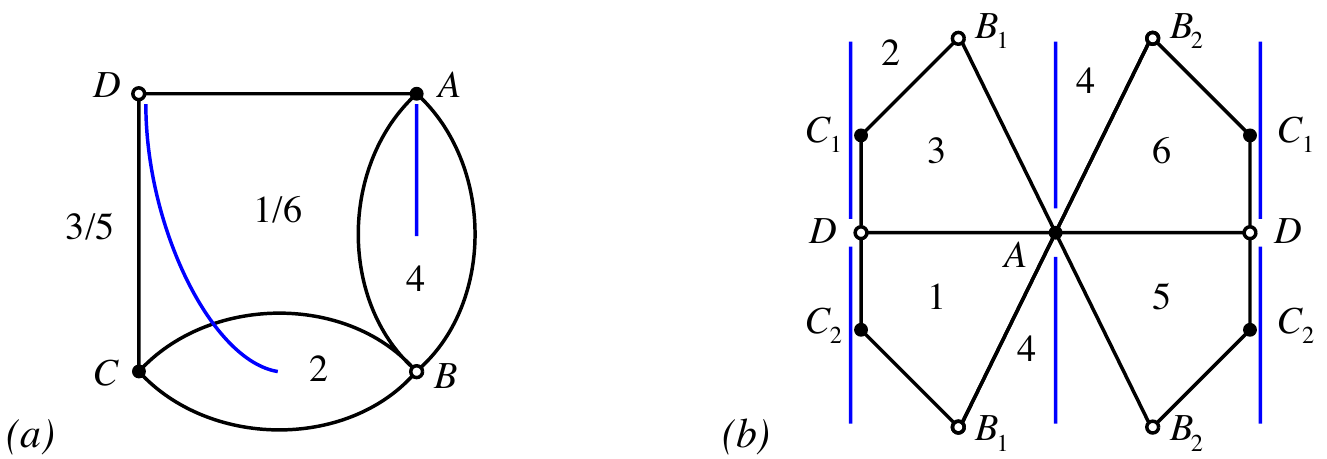}}
\end{picture}
\end{center}
\caption{The tiling \cite[(3.34)]{Davey:2009bp} from a genus 0 map}
\label{fig:dessin334}
\end{figure}

The genus 0 dessin is depicted in Figure \ref{fig:dessin334}\refpart{a},
together with two cuts between $A$, $D$ and centers of two cells (2 and 4).
Those two cells will have one pre-image on the genus 1 dessin, but with double valencies.
The real line is symmetrically represented by the straight rectangle $ABCD$. 
The $x$-coordinates of the vertices $A$, $B$, $C$, $D$ are, 
respectively, $x=2$, $x=0$, $x=-2$, $x=\infty$, and the branching centers 
of the cells $2$ and $4$ are, respectively, $x=1$ and $x=-1$.

The pulled-back dessin on a fundamental region of $E$ is depicted 
in Figure \ref{fig:dessin334}\refpart{b}. The $\IR$-symmetric cut from $D$ 
runs though the point $C$, hence the cuts are nearly parallel and
the genus 1 drawing is somewhat easier than in Figure \ref{fig:dessin328}\refpart{b}. 
One may consider the complementary branch cuts along the segment $A$, $D$ 
and between the cells 2, 4. This leads to considering Figure \refpart{b} in the vertical direction.

Figure \ref{fig:d334isog} depicts the doubly periodic dessin,
with the upper left gray square \refpart{a} representing a fundamental domain.
The dessin matches the tiling \cite[(3.34)]{Davey:2009bp}
after an interchange of the black and white vertices and a rotation by $\pi/2$.
The catalogue \cite{Davey:2009bp} gives exactly two 
tilings with the passport $[3^2\,6/4^3/2^6]$, up to colouring of the vertices.
As we see, both tilings can be obtained as compositions with a degree 6 genus 0 dessin.
They can be distinguished by the fact the degree 6 vertex in  Figure \ref{fig:dessin328}
is surrounded by all 6 distinct cells. 

\comment{
Generally, genus 1 dessins can be identified in the catalogue \cite{Davey:2009bp} 
by comparing the quiver diagrams, which is the oriented incidence (multi)graph of the dessin cells. 
The quivers can be drawn  rather fast  directly from the genus 0 dessins, 
as in Figure \ref{fig:dessin328}\refpart{a} or Figure \ref{fig:dessin334}\refpart{a}, by following incidences of cells, vertices on the genus 0 dessin.
The quiver for (\ref{eq:phi28}) is depicted in Figure \ref{fig:dessin328}\refpart{c}.
}

\begin{figure} 
\begin{center}
\begin{picture}(440,424)
\put(-6,-12){\includegraphics[height=432pt]{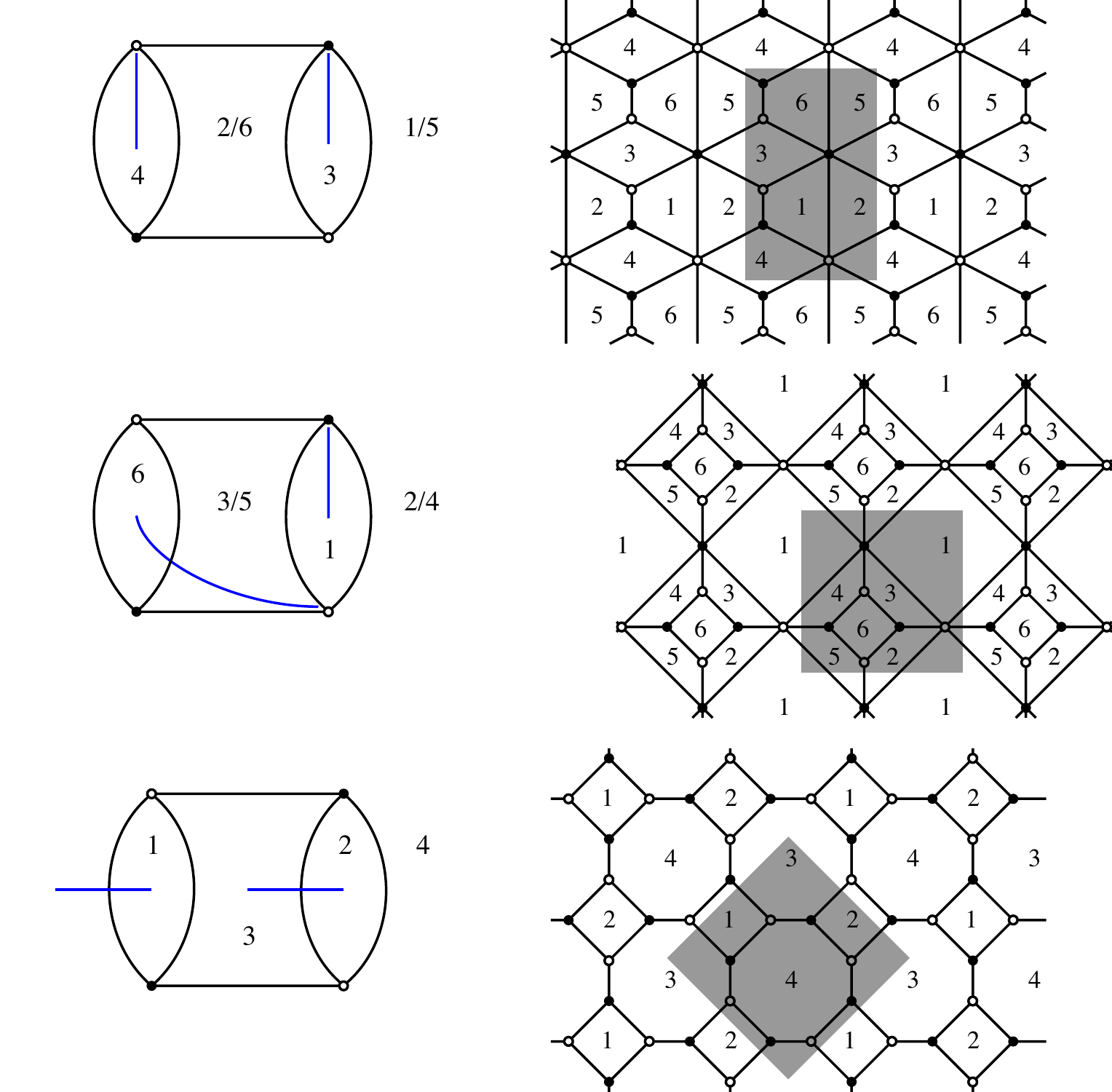}}
\put(12,290){$\refpart{a}$ Tiling \cite[\!(3.32)]{Davey:2009bp},
$\displaystyle j\!=\!\frac{52^3}3$.}
\put(12,144){$\refpart{b}$ Tiling \cite[(3.37)]{Davey:2009bp}, $j\!=\!1728$.}
\put(12,-2){$\refpart{c}$ Tiling $F_0$(II), $\!j=1728$.}
\end{picture}
\end{center}
\caption{Tilings from a genus 0 map with the passport $[3^2/3^2/2^21^2]$}
\label{fig:dessin332} 
\end{figure}

\subsection{A Multitude of New Degree 12 Maps}
\label{sec:multitude}

Encouraged by the ease with which we can give new maps by simple composition, 
consider the genus 0 Belyi map
\begin{equation} \label{eq:comp33}
\varphi_3=\frac{x^3\,(x-2)^3}{(2x-1)^3}=1+\frac{(x^2-4x+1)(x^2-x+1)^2}{(2x-1)^3}
\end{equation}
with the passport $[3^2/3^2/2^21^2]$ on $\IP^1$ and the dessin, together with the cuts are presented later on in Figure \ref{fig:dessin332}\refpart{a}.

On the elliptic curve
\begin{equation} \label{eq:ec52}
y^2=x\,(x^2-4x+1)
\end{equation}
with $j$-invariant $52^3/3$, the map $\varphi_3$ has passport  $[3^2\,6/3^2\,6/2^6]$ 
and the dessin is identified as \cite[(3.32)]{Davey:2009bp}, yet another new explicit Belyi pair.

On the other hand, had we chosen the elliptic curve
\begin{equation}\label{eq:ec1728-1}
y^2=(x-2)\,(x^2-4x+1)
\end{equation}
with $j$-invariant 1728,
the map $\varphi_3$ would have the same passport $[3^2\,6/3^2\,6/2^6]$ but the dessin would be different; it is in fact
\cite[(3.37)]{Davey:2009bp}, after following the cuts in Figure \ref{fig:dessin332}\refpart{b}.

Continuing under this light, on the genus 1 curve
\begin{equation}\label{eq:ec1728-2}
  y^2=(x^2-4x+1)(x^2-x+1) \ ,
\end{equation}
also with $j=1728$, $\varphi_3$ has the passport 
$[3^4/3^4/4^2\,2^2]$. 
The genus 0 and 1 dessins are 
given in Figure \ref{fig:dessin332}\refpart{c}.



\subsection{Genus 1 Maps of Degree 10}
\label{sec:deg10}

We now move onto maps of degree other than 12 as we encounter more complicated algebraic number.
Consider the genus 0 Belyi map
\begin{equation} \label{eq:deg10a}
\varphi_4=-\frac{x^2\,(4x+5)^3}{(5x+4)^3}=1-\frac{(x+1)\,(8x^2+11x+8)^2}{(5x+4)^3}
\end{equation}
with the passport $[3\;2/3\;2/2^2\,1]$. 
Let the genus 1 curve be
\begin{equation} \label{eq:ellc10a}
y^2=x\,(x+1)\,(x-\alpha) \ , \qquad\mbox{with}\quad
\alpha=-\frac{11+3\sqrt{-15}}{16}.
\end{equation}
Note that $\alpha$ is a root of $8x^2+11x+8$.
The expression $\varphi_4$ defines a degree 10 Belyi function on this curve,
with the passport $[3^2\,4/3^2\,4/2^3\,4]$.
The genus 1 dessin is identified as \cite[(3.12)]{Davey:2009bp}.
The genus 0 dessin (with branch cuts marked) and the genus 1 dessin 
are presented in Appendix figure \refpart{A}.  

The $j$-invariant of (\ref{eq:ellc10a}) equals $-108/5$, 
hence this curve can be defined over $\IQ$.
After the scaling transformation
\begin{equation} 
x\mapsto \frac{9+\sqrt{-15}}{48}\,x, \qquad
y\mapsto \frac{-7+\sqrt{-15}}{192}\, \textstyle \sqrt{-3+\sqrt{-15}}\;y
\end{equation}
we indeed obtain the elliptic curve $y^2=x\,(x^2+9x+24)$. 
But then the expression for the Belyi map $\varphi_4$ contains $\sqrt{-15}$ after this transformation.

A similar example is 
\begin{equation} \label{eq:comp5}
\varphi_5=\frac{x^2\,(x+5)^3}{(5x+1)^3}=1+\frac{(x-1)^3\,(x^2+18x+1)}{(5x+1)^3}
\end{equation}
which has the passport $[3^2\,4/3^2\,4/3^2\,2^2]$ on the curve
  $y^2=x\,(x^2+18x+1)$.
The genus 1 dessin is given in  Appendix figure \refpart{B},
together with the genus 0 dessin, the cut-data, the $j$-invariant,
and the identification with \cite[(3.5)]{Davey:2009bp}.

\section{Higher Degree Maps}\setall
\label{sec:higherd}

Continuing in the above fashion, we can move to higher degree maps.
In this section, we continue the above algorithm of establishing a genus 0 map $\varphi$, and finding an appropriate elliptic curve $E$ which covers the Riemann sphere and pulling $\varphi$ back to $E$.
A {\sf Maple} programme 
used in \cite{vHVHeun} was employed efficiently to find the maps 
presented below.

\subsection{Maps of Degree 14}
The genus 1 passport $[4^2\,3^2/4^2\,3^2/2^4\,3^2]$ can be obtained from
the following genus 0 passport by the following compositions:
\begin{itemize}
\item[\refpart{C}] $[2^2\,3/2^2\,3/2^2\,3]$; the Belyi map is
\begin{equation}\label{eq:C}
-\frac{(x+3)^3\,(x^2-x+2)^2}{(x-3)^3\,(x^2+x+2)^2}=1-\frac{2x^3\,(x^2+7)^2}{(x-3)^3\,(x^2+x+2)^2}.
\end{equation}
The genus 1 curve is $y^2=(x^2+7)(x^2+x+2)$ or  $y^2=(x^2-x+2)(x^2+x+2)$, 
isomorphic to $y^2=x(x+1)(x-7)$. 

\item[\refpart{D}] $[4\,3/4\,3/1^4\,3]$; the Belyi map is
\begin{equation}\label{eq:D}
-\frac{x^4\,(x-7)^3}{(7x-1)^3}=1-\frac{(x-1)^3\,(x^4-18x^3+90x^2-18x+1)}{(7x-1)^3}.
\end{equation}
The genus 1 curve is $y^2=x^4-18x^3+90x^2-18x+1$, 
isomorphic to $y^2=x(x^2+42x-7)$.

\item[\refpart{E}] $[4\,3/2^2\,3/1^2\,2\,3]$,  or equivalently, $[2^2\,3/4\,3/1^2\,2\,3]$.
These Belyi maps are defined over the cubic field 
$\IQ(\xi)$ with \mbox{$\xi^3-\xi^2+5\xi+15 = 0$}. Explicitly, a Belyi map is
\begin{equation}\label{eq:E}
\frac{x^3\,\big(x^2+14(1-\xi)x-7\xi^2-70\xi-245\big)^2}
{(14\xi^2-180\xi+90)\,\big(7x-4\xi^2+4\xi-8\big)^3}.
\end{equation}
The genus 1 curve is
\[
y^2=\big(x^2+14(1-\xi)x-7\xi^2-70\xi-245\big)
\big(x^2+(\xi^2-18\xi+13)x-8\xi^2+8\xi-144\big),
\] 
isomorphic to $y^2=x\,\big(x^2+(44\xi^2-125\xi+415)x-6\xi^2+4\xi-50\big)$.
\end{itemize}

The labels \refpart{C}--\refpart{E} refer to Appendix figures. 
In particular, two dessins \refpart{E1}, \refpart{E2} are given for
the cubic Galois orbit. The dessin \refpart{E2}
is defined over $\IR$ and is associated with the so-called $L^{333}$(II) theory.  


The genus 1 passport $[4^2\,3^2/4^2\,3^2/2^5\,4]$ can be obtained from
the following genus 0 passport by compositions:
\begin{itemize}
\item[\refpart{F}]  $[4\,3/4\,3/2^2\,1^3]$; the Belyi map is defined over $\IQ(\sqrt{7})$:
\begin{equation} \hspace{-10pt}\label{eq:F}
\frac{\big(x+4\sqrt7+7\big)^3\big(x-3\sqrt7+7\big)^4}
{\big(x-4\sqrt7-7\big)^3\big(x+3\sqrt7-7\big)^4}
=1+\frac{98\big(x^2\!+8\sqrt7-21\big)\big(x^2\!-2\sqrt7+28\big)^2}
{\big(x-4\sqrt7-7\big)^3\,\big(x+3\sqrt7-7\big)^4}.\!
\end{equation}
The elliptic curve is $y^2=(x-\beta)\big(x^2+8\sqrt7-21\big)$,
where $\beta=(\sqrt7+1)\sqrt{2\sqrt7-7}$ is a root of $x^2-2\sqrt7+28=0$.
Like with (\ref{eq:deg10a}), the map is defined over a quadratic extension of $\IQ(\sqrt7)$,
but either the elliptic curve or the map expression can be defined over $\IQ(\sqrt7)$.
The elliptic curve is isomorphic to 
\[
y^2=x\,\big(x^2+(10-2\sqrt7)x+69-24\sqrt7\big). 
\]
\item[\refpart{G}]  $[4\,3/2^2\,3/2^3\,1]$, or equivalently, $[2^2\,3/4\,3/2^3\,1]$.
A genus 0 Belyi map is
\begin{equation}\label{eq:G}
\frac{(x+9)^3\,\big(x^2-3x+4\big)^2}{x^3\,(x-7)^4}
=1+\frac{(7x^3-14x^2+63x-108)^2}{x^3\,(x-7)^4}.
\end{equation}
The elliptic curve is $y^2=(x-\lambda)(x^2-3x+4)$, 
where $\lambda$ is a root of $7x^3-14x^2+63x-108$.
This time the $j$-invariant is a cubic algebraic number.
Different cuts on a genus 0 dessin give different genus 1 dessins
like in \S \ref{sec:multitude};
see Appendix figures \refpart{G1} and \refpart{G2}. 
The cubic number field is more concisely defined by the polynomial
$\eta^3+\eta^2-2\eta+6$. This was found by computing the field discriminant ($-1176$)
and using the data base \cite{NFDB}.
\end{itemize}
The labels \refpart{F}--\refpart{G} refer to pairs of Appendix items. 

\subsection{Maps of Degree 16 and 18}
We can continue to higher degrees.
So far, we have seen degrees 10, 12, an 14,
in this subsection, we will see how degrees as high as 18 can be readily reached.

The genus 1 passport $[5^2\,3^2/4^4/2^8]$ of degree 16 comes from genus 0 passports 
$[5\,3/4^2/2^2\,1^4]$, $[5\,3/4\,2^2/2^3\,1^2]$ and $[5\,3/2^4/2^4]$.
There are no Belyi maps with the last passport by Remark \ref{rm:nonexist}.
The first passport gives the Belyi map
\begin{equation}\label{eq:16-1}
-\frac{5\,(x^2-10)^4}{4\,x^5\,(x-8)^3} = 1-\frac{(x^4-4x^3-8x^2+8x+20)\,(3x^2-10x+20)^2}{4\,x^5\,(x-8)^3}
\end{equation}
The applicable genus 1 curve is obvious; 
it isomorphic to $y^2=x^3-75x+290$ or $y^2=x^3+x^2-8x+8$.
The genus 0, 1 dessins are depicted in Appendix figure \refpart{H}.

The passport $[5\,3/4\,2^2/2^3\,1^2]$ has one Galois orbit
defined over $\IQ(\sqrt{10})$:
\begin{equation}\label{eq:16-2}
\frac{125\,x^4\,(x^2-4x-2\sqrt{10}-2)^2}{(7\sqrt{10}+25)\,(8x+2\sqrt{10}-3)^3}.
\end{equation}
The genus 1 curve is
\begin{equation}\label{ec:16-2}
y^2= \big( x^2-4x-2\sqrt{10}-2 \big) \,
\big( 5x^2-(6\sqrt{10}+10)x-4\sqrt{10}-5 \big),
\end{equation}
isomorphic to $y^2= x\,\big( x^2+(22\sqrt{10}+58)x-50\sqrt{10}-125 \big)$.
The dessins are depicted in  
Appendix figures \refpart{I1} and \refpart{I2}.

The genus 1 passport $[3^6/3^6/2^3\,3^2\,6]$ of degree 18 comes from the genus 0 passport 
$[3^3/3^3/1^3\,3^2]$. The Belyi map is
\begin{equation} \label{eq:deg18a}
\frac{(x^3+3x^2-3)^3}{27(x+1)^3(x+2)^3}=1+\frac{(x^3-9x-9)\,(x^2+3x+3)^3}{27(x+1)^3(x+2)^3}.
\end{equation}
The genus 1 curve is
$y^2=(x^3-9x-9)(x+1-\omega_3)$. This curve is isomorphic to $Y^2=X^3+1+\omega_3$ (hence $j=0$)
via 
\begin{equation}
x=  (\omega_3-1)\,\frac{X+1-\omega_3}{X+1}, \qquad 
y= \frac{3\omega_3 Y}{(X+1)^2}.
\end{equation}
The 
dessins are depicted in 
Appendix figure \refpart{J}. The genus 1 dessin can be redrawn with the hexagonal (a)symmetry.


\section{Isogenies}
\label{sec:isogenies}\setall
In the foregoing we have given a host of examples of obtaining genus 1 dessins by composing the easier genus 0 dessins
with a quadratic covering of the sphere.
A straightforward way to get a high degree genus 1 Belyi map is composition
of a genus 1 covering $\varphi_1:E_0\to \IP^1$
with a covering $\psi:E_1\to E_0$ of genus 1 curves.
The latter has to be unramified by the Hurwitz formula. 

If $E_0,E_1$ are given the group structure as elliptic curves,
and their group identity points are mapped to each other,
the map $\psi$ is an {\em isogeny}.
Computation of isogenies is a classical problem \cite{kohel}, \cite{isogn}, 
now especially actively explored for elliptic curves defined over finite fields 
in the light of applications to cryptography.
This section demonstrates comfortable interplay between
isogenies and special Belyi maps, suggesting a fresh computational approach
for both kind of functions.


As an example, 
we can compose the well-known degree 3 Belyi map \cite[(1.1)]{Davey:2009bp}
\begin{align} 
\label{eq:ecj0} 
\Phi_0 = \frac{1+Y}2 \qquad \mbox{on} \qquad Y^2 = X^3+1
\end{align}
with the degree 2 isogeny
\begin{equation}
X=\frac{(x-1)(x+3)}{4x}, \qquad Y=\frac{y\,(x^2+3)}{8x^2}
\end{equation}
from the $j=54000$ curve $y^2=x\,(x^2+6x-3)$, 
and get this degree 6 Belyi map:
\begin{equation}\label{deg6}
\frac12+\frac{y\,(x^2+3)}{16x^2}
\qquad \mbox{on} \qquad y^2=x\,(x^2+6x-3).
\end{equation}
The passport is $[3^2/3^2/3^2]$. 
The dessin is recognized as \cite[(2.1)]{Davey:2009bp}.

If a given genus 1 Belyi map on $Y^2=X^3+\ldots$ is a composition 
with genus 0 Belyi map  $\varphi_0(X)$ as considered in the previous section,
its composition with an isogeny from an elliptic curve in a Weierstra\ss{}
form $y^2=x^3+\ldots$ will have the same property, 
because the $X$-component of the isogeny is 
a rational function of $x$ as well known \cite[\S 2]{kohel}.
The Belyi maps (and the isogenies!) can be comfortably computed 
by considering them as compositions with genus 0 maps as above.
We demonstrate this on the isogenous versions for the map $\varphi_3$ in (\ref{eq:phi28})
and for the curves isogenous to $Y^2=X^3-X$.

\subsection{Degree 24 maps}
\label{sec:deg24}

\begin{figure}[p] 
\begin{center}
\begin{picture}(330,270)
\put(0,-3){\includegraphics[height=275pt]{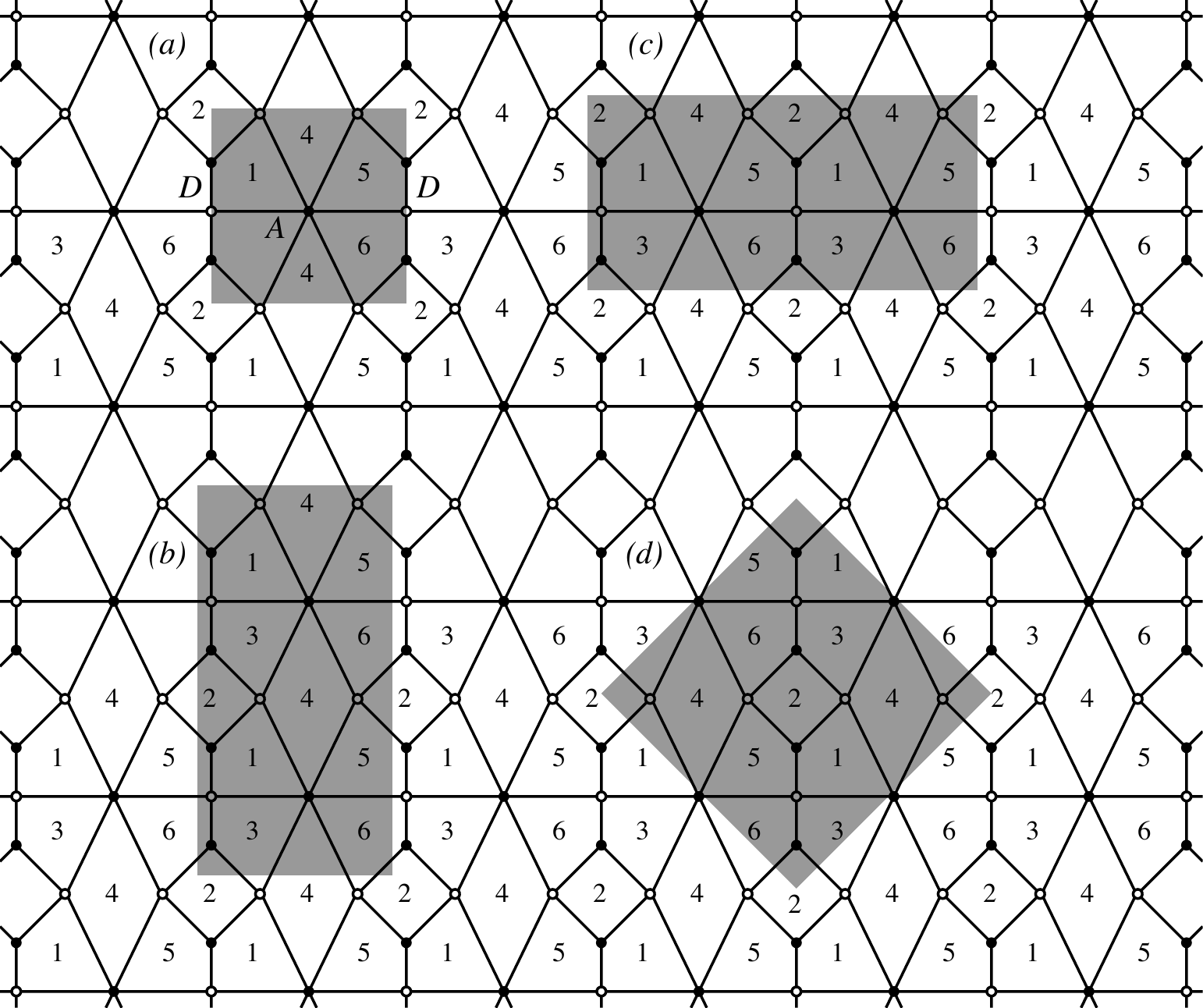}}
\end{picture}
\end{center}
\caption{The genus 1 dessins isogenous to \cite[(3.34)]{Davey:2009bp}}
\label{fig:d334isog}
\end{figure}

\begin{figure}
\begin{center}
\begin{picture}(440,274)
\put(-2,-3){\includegraphics[height=275pt]{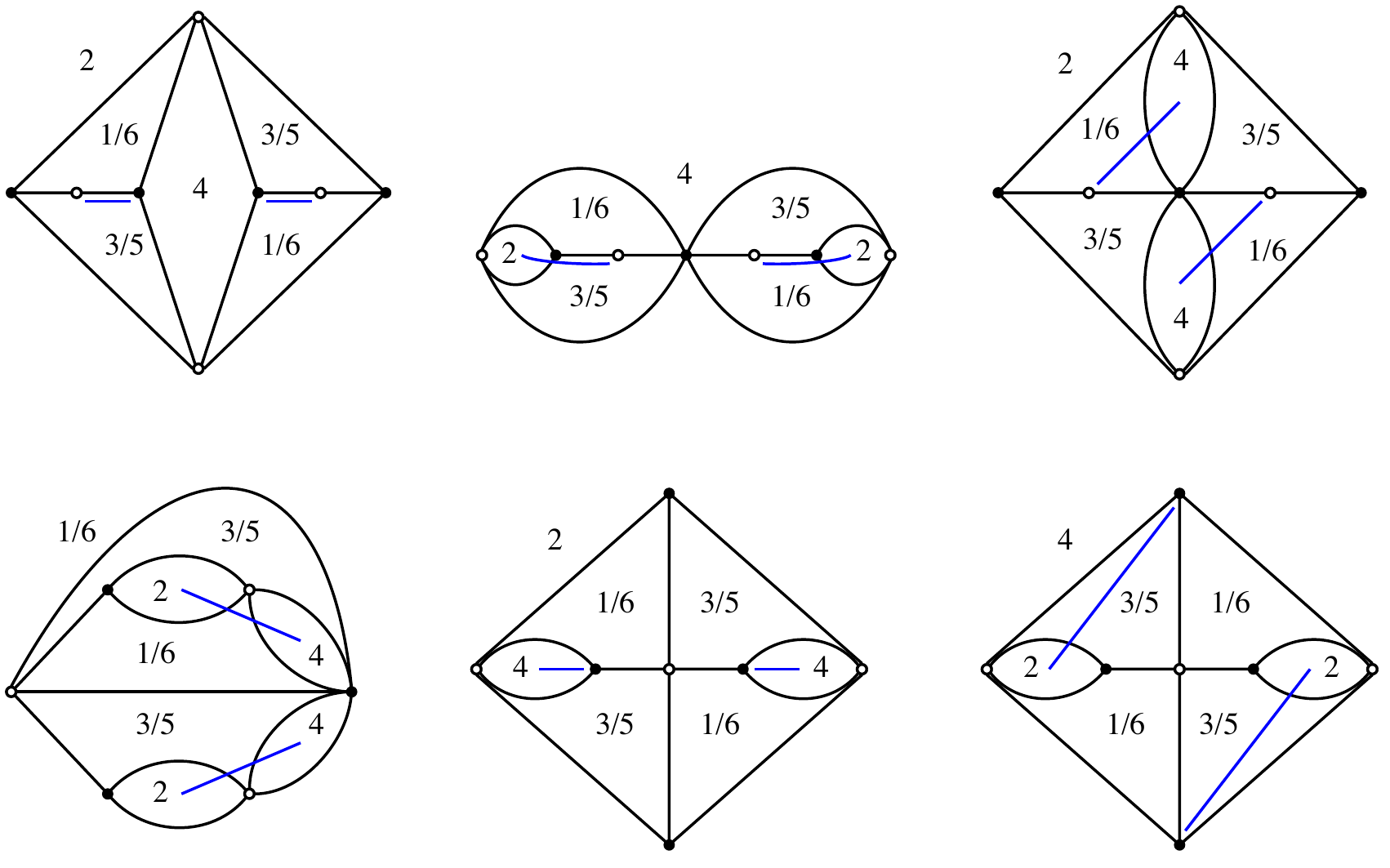}}
\put(0,148){\refpart{$\,b'$}}  \put(142,148){\refpart{$c'$}}   \put(309,148){\refpart{$d'$}} 
\put(0,-2){\refpart{$\,b''$}}  \put(142,-2){\refpart{$c''$}}   \put(309,-2){\refpart{$d''$}} 
\end{picture}
\end{center}
\caption{The genus 0 dessins corresponding to isogenies of \cite[(3.34)]{Davey:2009bp}}
\label{fig:d334isog0}
\end{figure}

As known from the theory of modular curves \cite{EMod},
isogenies of degree $N$ are generically counted by integer triples $(k,p,\ell)$
such that $kp=N$, $0\le \ell<k$.
In particular, an elliptic curve generically has three curves that are 2-isogenous to it.
There are three isogenies of degree 2 to the curve (\ref{eq:ec28}),
given by
\begin{equation} \label{eq:isog28}
x=\frac{X^2+4}{4X}, \qquad
x=\frac{X^2+6X+1}{4X} \quad\mbox{or}\quad
x=\frac{X^2+2X+9}{4X}.
\end{equation}
The three isogenous curves are, respectively,
\begin{equation} \label{eq:isog34}
Y^2=X\,(X^2-8X+4), \quad
Y^2=X\,(X^2+10X+1), \quad
Y^2=X\,(X^2-2X+9).
\end{equation}
Their $j$-invariants are $52^3/3,194^3/3,46^3/81$ respectively. 
The first curve happens to be isomorphic to  (\ref{eq:ec52}).
By making the substitutions (\ref{eq:isog28}) into (\ref{eq:phi28}) we get
three genus 1 Belyi maps with the passport $[3^4\,6^2/4^6/2^{12}]$.
The dessins look the same as doubly periodic tilings, 
but their fundamental domains have area  twice larger than Figure \ref{fig:d334isog} \refpart{a}.
The three fundamental domains are depicted in Figure \ref{fig:d334isog} \refpart{b}--\refpart{d}.

The substitutions (\ref{eq:isog28}) into (\ref{eq:phi28}) also give genus 0 Belyi maps
that determine the genus 1 maps by a quadratic covering.
The first genus 0 map has the passport $[3^4/4^2\,2^2/2^6]$,
while the other two have the passport $[3^2\,6/4^2\,2^2/2^5\,1^2]$. 
Their dessins are depicted in Figure 
\ref{fig:d334isog0}\refpart{$\,b'$},\refpart{$c'$},\refpart{$d'$}, respectively.
Each of them is obtained from two copies of 
Figure \ref{fig:dessin334}\refpart{a} glued along a single branch cut.
For example, Figure \ref{fig:d334isog0}\refpart{$\,b'$} is obtained 
after a branch cut between the cells 2, 4 in Figure \ref{fig:dessin334}\refpart{a}. 
The corresponding genus 1 dessin  is obtained by cutting  
Figure \ref{fig:d334isog0}\refpart{$\,b'$} along two pre-images of the segment $AD$.
It has the fundamental domain \refpart{b} in Figure \ref{fig:d334isog}. 
Similarly, the dessins \refpart{$c'$}, \refpart{$d'$} in Figure \ref{fig:d334isog0}  
are obtained by cutting Figure \ref{fig:dessin334}\refpart{a} between vertex $A$
and the cells 4 or 2 (respectively). The corresponding genus 1 dessins 
have the fundamental domains \refpart{c}, \refpart{d} in Figure \ref{fig:d334isog}.

The same three genus 1 dessins are obtained by cutting Figure \ref{fig:dessin334}\refpart{a}
first along the ``parallel" cuts from the vertex $D$ to (respectively) the vertex $A$
or the cells 2 or 4. The intermediate genus 0 dessins are depicted in 
Figure \ref{fig:d334isog0}\refpart{$\,b''$},\refpart{$c''$},\refpart{$d''$}, respectively.
The dessins \refpart{$c''$},\refpart{$d''$} are the same by the symmetry of Figure \ref{fig:dessin334}\refpart{a},
but their branch cuts to genus 1 dessin differ. 
The two different Belyi maps can be obtained by substituting $x=X^2+2$ or $x=X^2+1$
into  (\ref{eq:phi28}). The isogenous elliptic curves are then
\begin{equation}
Y^2=(X^2+1)(X^2+3), \quad
Y^2=(X^2+2)(X^2+3), \quad
Y^2=(X^2-1)(X^2+2),
\end{equation}
respectively isomorphic to (\ref{eq:isog34}). 
The apparent duality of Belyi maps in Figure \ref{fig:d334isog0} 
reflects commutativity of constructing Riemann surfaces by "parallel" cuts
on the Riemann sphere.

More generally, a degree 2 isogeny map to $y^2=x\,(x^2+a\,x+b)$ is given by $x=X^2$,
with the isogenous curve being $Y^2=X^4+a\,X^2+b$. 
This 2-descent goes back to Fermat's proof (cf.~\cite{Beu})
on integer solutions of $a^4+b^4=c^4$.
By (\ref{eq:ec4to3}), the isogenous curve is isomorphic to $Y^2=X\,(X^2-2a\,X+a^2-4b)$.

\subsection{Isogenies on the Square Lattice}
\label{sec:square}

\begin{figure} 
\begin{center}
\begin{picture}(440,490)
\put(-2,-4){\includegraphics[width=440pt]{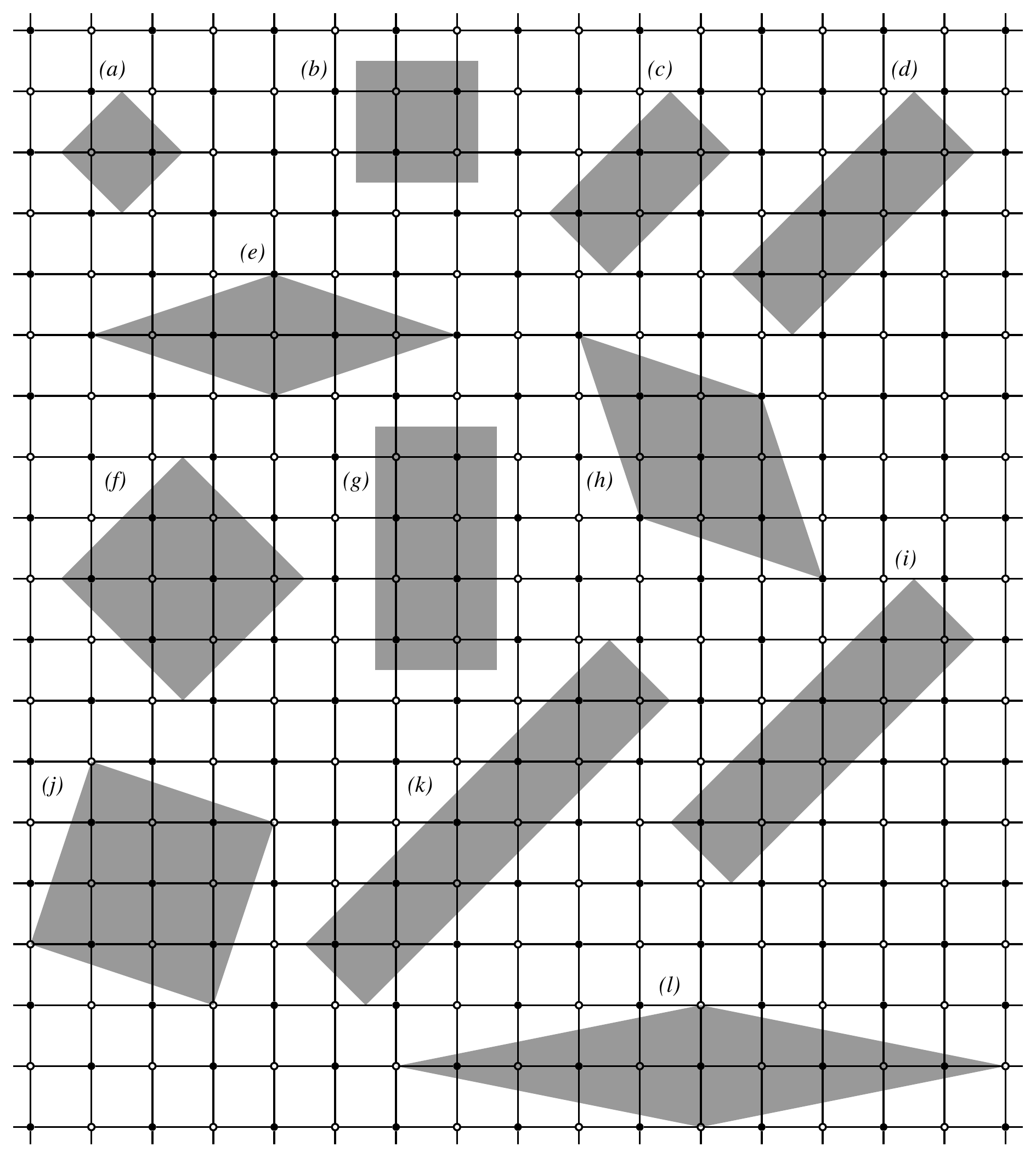}}
\end{picture} 
\end{center}
\caption{   
  \comment{The genus 0 dessins for \cite[(3.32), (3,37)]{Davey:2009bp} 
        and $\IF_0$(II), etc., 
        Given a square tiling of the doubly periodic plane , we can chose the fundamental domain (drawn as shaded region) to give different dessins.  The simplest is (a) where there is a single pair of black-white nodes in the fundamenal domain.}       
Compositions of isogenies to $Y^2=X^3-X$ and the map $\Psi_0=X^2$
have the simple square lattice (with alternating black/white nodes) 
as the doubly periodic dessins d'enfant. An isogeny of degree $k$ 
gives a composition of degree $4k$ and a sublattice of index $2k$,
represented by a (shaded) fundamental domain of the lattice area $2k$.
} \label{fig:SquareDessin}
\end{figure}


Let $E_1$ denote the elliptic curve $Y^2=X^3-X$ with $j=1728$.
This subsection considers compositions of isogenies to $E_1$ with the map 
$\Psi_0=X^2$. 

The map $\Psi_0$ is a degree 4 Bely map on $E_1$
with the branching passport $[4/4/2^2]$. Its dessin d'enfant forms a simplest
square tiling of the plane in Figure \ref{fig:SquareDessin}, 
made of entirely valency 4 alternating black/white nodes.
The fundamental domain in \refpart{a} contains a single pair of black and white nodes 
This dessin is identified as the conifold $L^{111}$ in \cite[(1.2)]{Davey:2009bp}.


A composition of an isogeny of degree $k$ with $\Psi_0$ 
is a genus 1 Belyi map of degree $4k$, with the 
passport $[4^k/4^k/2^{2k}]$.
As mentioned just before \S \ref{sec:deg24}, 
these maps are compositions of a genus 0 map
and a quadratic covering $E\to\IP^1$.
The dessins d'enfant make up the same square lattice,
but the area of their fundamental regions (parallelograms) is $k$ times larger.
The number of distinct isogenies of degree $k$ equals the number of distinct sublattices
of the square lattice of index $2k$.
We show various possibilities in Figure \ref{fig:SquareDessin}.
The shaded region shows a choice of the fundamental domain 
when drawing the dessin on the torus. 

There are two distinct isogenies of degree 2 onto $E_1$:
\begin{align} \label{eq:isog1a}
& X = \frac{x^2+1}{2x}, && Y = \frac{ y \, (x^2-1)}{2\sqrt2 \, x^2}
&& \mbox{from } y^2=x^3+x;
\\[2pt] \label{eq:isog1b}
& X = \frac{(x+1)^2}{4x}, &&
Y = \frac{y\,(x^2-1)}{8x^2} && \mbox{from } y^2=x\,(x^2+6x+1); 
\end{align}
these curves have $j=1728$ and $66^3$ respectively and yield, 
via composition with $\Psi_0$, Belyi maps of degree 8:
\begin{align} \label{psi11}
\psi_1 = \frac{(x^2+1)^2}{4x^2} = & \, 1+\frac{(x^2-1)^2}{4x^2} 
&& \mbox{on } y^2=x^3+x; \\[2pt]
\label{psi12}
\psi_2 = \, \frac{(x+1)^4}{16x^2} \, =  & \, 1+\frac{y^2(x-1)^2}{16x^3}  
&& \mbox{on } y^2=x\,(x^2+6x+1) \ .
\end{align}
As genus 0 Belyi maps of degree 4, they have the passports $[2^2/2^2/2^2]$
and $[2^2/4/1^2\,2]$, respectively. 
The map $\psi_1$ can also be obtained from a genus 0 map with the passport $[4/4/1^4]$,
as $\psi_1=X^4$ on the curve $Y^2=X^4-1$.
The matching isomorphism of the elliptic curves 
is given by $X=y/(\sqrt{2}\,x)$. The elliptic curve in (\ref{eq:isog1a}) can be rescaled 
$x\mapsto -i x$, $y\to\sqrt{i}\,y$ to be identical to $E_1$; the isogeny is then an endomorphism of $E_1$.
The rescaling interchanges the expressions for $\psi_1$ and $1-\psi_1$.

The fundamental regions for $\psi_1,\psi_2$ are depicted 
in Figure \ref{fig:SquareDessin}\refpart{b}, \refpart{c}, respectively.
We may write $\psi_2=(X^2-1)^2$ with $X=y/2x$, using the simplest
expression for the genus 0 Belyi map with the relevant passport is $[4/2^2/1^2\,2]$.
The genus 1 curve is then $Y^2=(X^2-1)(X^2-2)$. 

To find 3-isogenies to $E_1$, we 
look for genus 0 Belyi maps
with the passports $[4\,2/4\,2/2^2\,1^2]$ or $[4\,2/2^3/2^3]$,
so that the genus 1 passport $[4^3/4^3/2^6]$ could be obtained
by a quadratic covering with 4 branching points. 
There are no Belyi maps with the second passport by Remark \ref{rm:nonexist}.
The first passport gives a Galois orbit defined over $\IQ(\sqrt3)$.
As a genus 1 map we get
\begin{equation} \label{eq:isog1c}
\psi_3=\frac{\sqrt3}{2} \frac{x^2(x+3\sqrt3-3)^4}{(3x+5\sqrt3-9)^4} 
\quad \mbox{on} \quad
y^2=x \, \big(x^2+8\sqrt3\,x-14\sqrt3+24\big).
\end{equation}
To get an isogeny without additional algebraic extensions, we rescale $E_1$
to $Y^2=X^3-6\sqrt{3}\,X$. The isogeny from the elliptic curve in (\ref{eq:isog1c}) is 
then
\begin{equation*}
X=\frac{3x(x+3\sqrt3-3)^2}{(3x+5\sqrt3-9)^2}, \quad
Y= \frac{3\sqrt{3}\,y \big(x+3\sqrt3-3\big) \big(x^2+(2\sqrt3-6)\,x-14\sqrt3+24\big)}{(3x+5\sqrt3-9)^3}.
\end{equation*}
The elliptic curve has two real components,
thus a fundamention region  for $\psi_3$ is a rectangle; 
it depicted in Figure \ref{fig:SquareDessin}\refpart{d}. 
The conjugated version $\sqrt3\mapsto -\sqrt3$ gives other Belyi map.
The conjugated elliptic curve has just one real component,
and a fundamental region is a rhombus as in Figure \ref{fig:SquareDessin}\refpart{e}.

The 4-isogenies to $E_1$ can be similarly found from genus 0 Belyi maps
with the passports $[4^2/4^2/2^2\,1^4]$, $[4^2/4\,2^2/2^3\,1^2]$, 
$[4\,2^2/4\,2^2/2^4]$, $[4^2/2^4/2^4]$ (up to permutation of the $0$ and $\infty$ fibers),
aiming for the genus 1 passport $[4^4/4^4/2^8]$.
As in \S \ref{sec:deg24}, 
some of the genus 0 Belyi maps with different passports lead to the same genus 1 map.
Four different genus 1 Belyi maps (or isogenies) are found:
\begin{enumerate}
\item[\refpart{f}] The map $\psi_1(x)^2$ on $y^2=x^4+6x^2+1$, \ or $\psi_1(x^2)$ on $y^2=x^4+1$. 
Either of these maps corresponds to the multiplication-by-2 endomorphism on $E_1$:
\begin{equation}\label{4iso}
(X,Y)\mapsto \left( \frac{(X^2+1)^2}{4\,(X^3-X)}, \;
\frac{(X^2+1)(X^4-6X^2+1)}{8\,Y^3} \right),
\end{equation} 
via $X=\big(x^2+1+\sqrt{x^4+6x^2+1}\big)/2$ or $X=x^2+\sqrt{x^4+1}$, respectively. 
\item[\refpart{g}] 
The map $\displaystyle -\frac{(x^2+6x+1)^2\,(x-1)^4}{64\,x^2\,(x+1)^4}$ on $y^2=x(x^2+6x+1)$, \\[2pt]
or $\displaystyle -\frac{(x^2+2x-1)^4}{64x^4}$ 
on $y^2=x^4+8x^3+18x^2-8x+1$.\\[2pt]
These maps are related by  $x\mapsto \big(-2x+\sqrt{x(x^2+6x+1)}\big)/(x+1)$. \\
The isogeny is a composition of (\ref{eq:isog1a}) and the twist $X\mapsto iX$ of (\ref{eq:isog1b}).
\item[\refpart{h}] The map 
$\displaystyle \frac{(17+12\sqrt{2})\,\big(x^2+(10-8\sqrt2)x+1\big)^4}{256\,x^2\,(x-1)^4}$ 
defined over $\IQ(\sqrt2)$, \\[2pt]
on the curve \mbox{$y^2=x \big(x^2+(66-48\sqrt2)x+1\big)$.} 
\item[\refpart{i}] The conjugated $\sqrt2\to -\sqrt2$ version of the previous map.
\end{enumerate}
The labels \refpart{f}--\refpart{i} refer to the fundamental regions in Figure \ref{fig:SquareDessin}.

We finish with the Belyi maps corresponding to isogenies of degree 5. 
The applicable genus 0 passports are $[4^2\,2/4^2\,2/2^4\,1^2]$ and $[4^2\,2/4\,2^3/2^5]$.
There are no Belyi maps for the latter passport by Remark \ref{rm:nonexist}.
The first one gives two Galois orbits,
defined over $\IQ(\sqrt{5})$ and $\IQ(i)$. A Belyi map defined with $\sqrt5$ is
\begin{equation}\label{ec:kl}
\frac{(9-4\sqrt5)\,x^2\,\big(x^2-(5\sqrt5+25)x-10\sqrt5+25\big)^4}
{\sqrt5\,\big(5x^2+(15\sqrt5-25)x-38\sqrt5+85\big)^4}
\end{equation}
defined on $y^2=x\big(x^2-(48\sqrt5+120)x-9\sqrt5+20\big)$.
A fundamental region is depicted in Figure \ref{fig:SquareDessin}\refpart{k}.
The conjugate $\sqrt5\to-\sqrt5$ version gives the fundamental region 
in Figure \ref{fig:SquareDessin}\refpart{l}.

The Belyi maps defined over $\IQ(i)$ are defined on $E_1$ and 
come from the endomorphism isogeny
\begin{equation} \label{eq:isogm2}
(X,Y) \mapsto \left( \frac{X\,(X^2-1-2i)^2}{((1+2i)X^2-1)^2}, 
\frac{Y\,(X^2-1-2i)\,(X^4+(2+8i)X^2+1)}{((1+2i)X^2-1)^3} \right).
\end{equation}
on $E_1$, and the complex conjugation. 
A corresponding fundamental region is depicted in Figure \ref{fig:SquareDessin}\refpart{j}.
Its mirror image is the conjugated Belyi map.
Recall that $E_1$ has {\em complex multiplication} \cite[\S 3]{kohel} by $i$,
meaning that the isogenies $E_1\to E_1$ form a ring 
isomorphic to 
$\IZ[i]$. The degree of such isogeny equals the
squared norm of the corresponding Gaussian integer.
These isogenies (and the corresponding Belyi maps) can be computed easily
by using the addition law on $E_1$. 
As described in \cite[\S 8]{vhpg}, \cite{vhpgell},
the isogenies also give transformations of the elliptic integral
\begin{equation}\label{F21}
\hpg{2}{1}{\frac12,\,\frac14}{\frac54}{z} = \frac{z^{-\frac14}}{2}\, 
\int_{\frac{1}{\sqrt{z}}}^{\infty}\,\frac{dX}{\sqrt{X^3-X}}.
\end{equation}
For example, 
\begin{equation} \label{eq:hpg5}
\hpg{2}{1}{1/2,\,1/4}{5/4}{z} = \frac{1-z/(1\!+\!2i)}{1-(1\!+\!2i)z}\;
\hpg{2}{1}{1/2,\,1/4}{5/4}{\frac{z\,(z-1-2i)^4}{\big((1\!+\!2i)z-1\big)^4}}.
\end{equation}

\section{Compositions with Other Genus 1 Coverings}
\label{sec:othercomp}
\setall
In the previous sections, we have considered composing genus 0 Belyi maps with unbranched covering of $\IP^1$ by an elliptic curve, as well as obtaining new genus 1 maps by composing with isogenies of an elliptic curve.
Now, the $j=0$ elliptic curve has very special symmetries and we wish to harness them.
In this seciton, we consider a parametric covering of degree 3, and its special case of maps to $j=0$.
These will produce yet more new classes of genus 1 maps.


\subsection{A Parametric Degree 3 Covering}

With $u$ as a parameter, consider a function
\begin{equation}
\Psi_1=Y+X+u \qquad \mbox{on} \qquad Y^2=X^3+(X+u)^2.
\end{equation}
This is a map of degree 3, because a generic equation $\Psi_1=v$ 
defines a line in $\IP^2$ that intersects the cubic curve at 3 points. 
Here are the branching points of $\Psi_1$:
\begin{enumerate}
\item The point at infinity, of branching order 3, where $\Psi_1=\infty$;
\item The point $(X,Y)=(0,-u)$ of branching order 3, where $\Psi_1=0$;
because this is a flex point on the cubic curve, 
and the line $\Psi_1=0$ is the tangent to it;
\item Two points $(X,Y)=(\frac34v-\frac32u,\frac14v+\frac12u)$
in the fibers $\Psi_1=v$ satisfying 
\begin{equation} \textstyle
(v-2u)^2+\frac{32}{27}\,v=0.
\end{equation}
\end{enumerate}
We modify $\Psi_1$ so that the branching fibers in \refpart{iii} would be roots
of any prescribed polynomial $v^2+Av+B$. This is achieved by defining
\begin{equation} \label{eq:cpsi2}
\Psi_2=\frac{Y+3X+\sqrt{B}}{2} 
\qquad \mbox{on} \qquad 
Y^2=\frac{16\,X^3}{A+2\sqrt{B}}+\big(3X+\sqrt{B}\big)^2. 
\end{equation}
The critical values $\Phi_2$ are $\infty$, $0$ and the roots of $v^2+Av+B=0$.
Note that $A+2\sqrt{B}=0$ only if the discriminant of $v^2+Av+B$ equals zero.

Having obtained the specialized form of $\Psi_2$, we now demonstrate its composition
with the following genus 0 Belyi map to generate several more new genus 1 dessins: 
\begin{equation} \label{eq:phi422}
\frac{z\,(z+4)^3}{4\,(2z-1)^3}=1+\frac{(z^2-10z-2)^2}{4\,(2z-1)^3},
\end{equation}
with passport $[3\,1/3\,1/2^2]$.
Keeping $\infty,0$ as the branching points in \refpart{i}--\refpart{ii},
we pick two branching fibers $v_1$, $v_2$ in \refpart{iii} 
among the roots of $(v+4)(2v-1)(v^2-10v-2)=0$
and make several 
choices:
\begin{itemize}
\item[\refpart{K}]  The choice of $v_1=-4$, $v_2=1/2$ gives $A=7/2$, $B=-2$.
An alternative compact expression for the elliptic curve in  (\ref{eq:cpsi2})
is obtained after the substitution $X=-(3x+1+\sqrt{-2})/2$, $Y=(1+2\sqrt{-2})y/2$  to arrive at
\begin{equation} \label{eq:ecK}
3y^2=16x^3+(4\sqrt{-2}-5)x^2-(2\sqrt{-2}+14)x-2\sqrt{-2}-3.
\end{equation}
Accordingly, substituting $z=\big((1+2\sqrt{-2})y-9x-\sqrt{-2}-3\big)/4$ into (\ref{eq:phi422})
we get a 
genus 1 Belyi function with the passport $[3^2\,6/3^2\,6/2^6]$.
This passport has already appeared twice in \S \ref{sec:multitude}, but here, 
we see that the dessin is identified with \cite[(3.35)]{Davey:2009bp}.

\item[\refpart{L}] 
The choice $A=-10$, $B=-2$ gives a Belyi map with the passport $[3^4/3^4/4^2\,2^2]$,
  as with the last case of \S \ref{sec:multitude}.
  After the substitution $X=(1+\sqrt{-2})(2x-\sqrt{-2})$, $Y=2y$
  into (\ref{eq:cpsi2}) we get the elliptic curve
  \begin{equation}  \label{eq:ecL}
  y^2=16x^3-(6\sqrt{-2}+9)x^2+(12\sqrt{-2}+6)x-2\sqrt{-2}+7 \ .
  \end{equation}
  Substituting $z=y+(3\sqrt{-2}+3)x-\sqrt{-2}+3\big)/4$ into (\ref{eq:phi422})
  gives the Belyi map on the above elliptic curve.

\item[\refpart{M}]  The other choices (such as $v_1=\!-4$, $v_2=5-3\sqrt3$)
give Belyi maps with the passport $[3^2\,6/3^4/2^4\,4]$. 
They form Galois orbits defined with field extension by  
$\sqrt{3\sqrt3-5}=(2-\sqrt3)\sqrt{\sqrt3+1}$.
\end{itemize}


\begin{figure}
\begin{picture}(440,158)
\put(-1,20){\includegraphics[width=440pt]{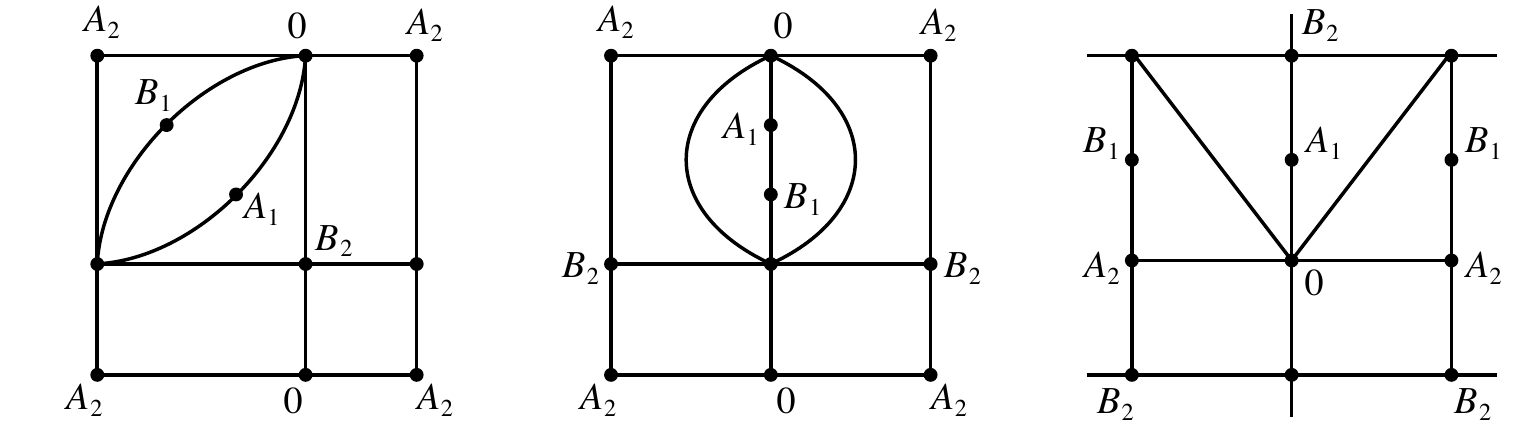}}
 \put(11,67){$\infty$} \put(123,67){$\infty$}
 \put(223,55){$\infty$} 
 \put(319,132){$\infty$}  \put(411,132){$\infty$}  \put(374,23){$\infty$} 
\put(0,2){\refpart{a}}  \put(145,2){\refpart{b}}   \put(292,2){\refpart{c}} 
\end{picture}
\caption{Inverse images of $\IR$ under the cubic map $\Psi_1$ with real $v$-fibers}
\label{fig:CubicInv}
\end{figure}

The dessins for the above are drawn in the Appendix, with the same labels
\refpart{K}--\refpart{M}.
To draw the dessins following the composition, we need to figure out the inverse image of the real line (with four marked points) under $\Psi_1$. 
The inverse image would be a graph (with the vertices marked by four ``colours")
on a torus with the following properties:
\begin{itemize}
\item It has six vertices: two vertices of order 6 (the preimages of $x=0$ and $x=\infty$), two vertices $A_2$, $B_2$ of order 4 (the other branching points), 
and two vertices $A_1$, $B_1$ of order 2 (the other two points in the branching fibers);
\item The edges divide the torus into 6 cells, 
each representing (as a preimage) a  half $\IC$-plane. 
Each cell has four vertices, of all four special fibers (``colours");
\item The colours of vertices with a common edge have a homotopical
``Schwarz reflection" symmetry. 
\end{itemize}

In Figure \ref{fig:CubicInv}, we present the possibilities of partitioning the fundamental domain of a torus into 6 cells of degree 4;
they give us all dessins for the maps in \refpart{K}--\refpart{M},
so the list is apparently complete for the cases of real roots of $v^2+Av+B$.
Dessin \refpart{a} applies to the cases when $x=0$, $x=\infty$ are not adjacent on $\IR$ in the cyclical sequence of the four points (thus $B<0$),  
while the dessins \refpart{b} and \refpart{c} apply to the adjacent case (thus $B>0$).
The dessin \refpart{c} is not a fundamental domain for the dessin; 
to have the doubly-periodic tiling, the square in Figure \refpart{c} 
has to be fitted in a brick fashion. The actual fundamental domain
gives awkward cell identifications across its edges, as we shall see in the Appendix.

Thus prepared, we can draw the dessins for the above 3 maps.
The figures are given in the Appendix figures and are made by fitting/overlaying
a half $\IC$-plane of the degree 4 dessin for (\ref{eq:phi422}) 
into each cell of a relevant dessin in Figure \ref{fig:CubicInv}. 
This amounts to fitting a shape in the form of letter ``N'' into each cell, 
with the end points at the $3$-branching nodes above $x=0$, $x=\infty$.
This is because the degree 4 (genus 0) dessin
is cut at $x=0$, $x=\infty$  
and at various other two points on $\IR$. 
In particular, for the dessin in \refpart{K} in the Appendix, the bending points of the letter ``N'' are placed at $A_i,B_i$ ($i\in\{1,2\}$) in each cell.
In dessin \refpart{L}, the ``N" bends at the middle points 
of $0A_i$ and $\infty B_i$. The dessins \refpart{M1}--\refpart{M3} 
follow a mixed case of bending at $A_i$ and the middle points $\infty B_i$ or $\infty0$. 

\subsection{Maps to $j=0$}
\label{sec:mapsj0}

Many examples of dessins with $j=0$ are obtained by composing genus 0 map 
\begin{equation}
  \Phi_0=z \quad \mbox{ on } \quad y^3=z(z-1) \ ,
\end{equation}
with the passport $[3/3/3]$. 
This map is equivalent to $\Phi_0$ in (\ref{eq:ecj0})
and is a well-known example in the physics literature \cite{Hanany:2011ra,Jejjala:2010vb}.

\begin{figure}
\begin{picture}(440,152)
\put(-3,2){\includegraphics[width=440pt]{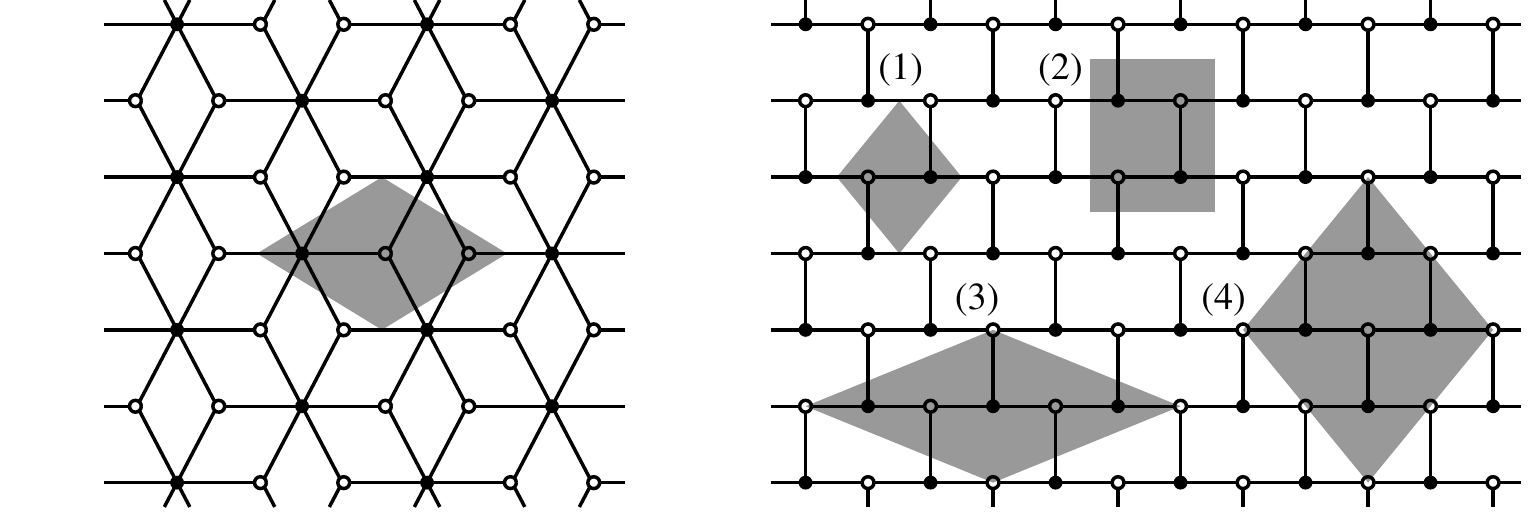}}
\put(2,8){\refpart{a}}  \put(194,8){\refpart{b}} 
\end{picture}
\caption{Dessins on the $j=0$ or $j=54000$ curves}
\label{fig:HexDessin}
\end{figure}

Considering $\varphi_0=z^2$ 
on the curve $y^3=(z-1)(z+1)$ gives a genus 1 map with the passport $[6/2^3/3^2]$. 
Since $\varphi_0-1=y^3$, 
it is also a composition of a degree 3 genus 0 Belyi function 
with the quadratic covering $z^2=y^3+1$. The dessin of $1-\varphi_0$ 
is depicted in Figure  \ref{fig:HexDessin}\refpart{a}.

Consider now the degree 5 map 
\begin{equation} \label{eq:deg15}
\frac{z\,(z-1-2i)^4}{\big((1\!+\!2i)z-1\big)^4}=1+\frac{(z-1)(z^2+(2+8i)z+1)^2}{\big((1\!+\!2i)z-1\big)^4}
\end{equation}
related to isogenies on the $j=1728$ curves; see (\ref{eq:hpg5}) and (\ref{eq:isogm2}).
On $y^3=z(z-1)$ this is a degree 15 map with the passport $[4^3\,3/4^3\,3/2^6\,3]$.
The genus 0 and 1 dessins are shown in Appendix figure \refpart{N}, respectively 
on the left and on the right sides.  To get the genus 1 dessin, we cut the $\IC$-plane
of the genus 0 dessin into two ``halfs" by drawing a circle through the three 
branching points (of different colours). 
The preimage of a ``half-plane" will be a curvilinear triangle (with the angles $\pi/3$).
We represent preimages of both ``half-planes" by two equilateral triangles 
with respective parts of the genus 0 dessin homotopically drawn of them. 
The genus 1 dessin is obtained by fitting the marked triangles into the triangular lattice,  
by alternating the two kinds of triangles in a vertex-matchin manner.

Finally, consider the degree 18 map
\begin{equation} \label{eq:deg18b}
\frac{z^3(z-2)^3}{(2z-1)^3}=1+\frac{(z^2-4z+1)(z^2-z+1)^2}{(2z-1)^3}
\end{equation}
on $y^3=(z^2-4z+1)\big(z+\omega_3\big)$. 
It has the passport $[3^6/3^6/2^3\,3^2\,6]$.
This is the same map as (\ref{eq:deg18a}),
via this isomorphism to $Y^2=X^3+1+\omega_3$:
\begin{equation}
z=\frac{\omega_3\,(1-Y)}{Y+\omega_3+1}, \quad 
y=\frac{(1-\omega_3)X}{Y+\omega_3+1}.
\end{equation}

\subsection{Isogenies on $j=0$}
\label{sec:isogj0}

The map $\varphi_3=z^3$ on the curve $y^3=z^3-1$ gives a degree 9 map
with the passport $[3^3/3^3/3^3]$.  After applying the isomorphism 
\[
X=-\frac{2^{2/3}\,y}{z-1}, \qquad Y=\sqrt{-3}\, \frac{z+1}{z-1},
\]
we get the 3-isogeny endomorphism on $Y^2=X^3+1$:
\begin{equation}
(X,Y)\mapsto \left( -\frac{X^3+4}{3X^2}, \; \frac{Y\,(X^3-8)}{3\sqrt{-3}\,X^3} \right).
\end{equation}
Alternatively, we can match $\Phi_0=(1+Y)/2$ with  
\[ 
\frac{\varphi_3}{\varphi_3-1}
=\frac{\big(Y+\sqrt{-3}\big)^3}{6\sqrt{-3}\,\big(Y^2-1\big)},
\]
and get the same isogeny in the form
\begin{equation}
(X,Y)\mapsto   \left( -\frac{Y^2+3}{3X^2}, \; \frac{Y\,(Y^2-9)}{3\sqrt{-3}\,(Y^2-1)} \right).
\end{equation}
Similarly, a genus 0 map with the passport $[3\,1/3\,1/3\,1]$
directly gives the multiplication-by-2 isogeny
\begin{equation} \label{eq:4isog}
(z,y)\mapsto \left( -\frac{z(z-2)^3}{(2z-1)^3}, \, 
\frac{y\,(z-2)(z+1)}{(2z-1)^2}\right)
\end{equation}
on $y^3=z\,(z-1)$. The translation onto $Y^2=X^3+1$ is straightforward.

Similarly as in Figure \ref{fig:SquareDessin}, the isogenous maps to $\Phi_0$ 
can be represented by sublattices of the hexagonal lattice.
Figure \ref{fig:HexDessin}\refpart{b} rather depicts the equivalent brick lattice,
with four fundamental regions representing, respectively, 
the maps $\Phi_0$, \eqref{deg6}, $\varphi_3$ and the isogeny \eqref{eq:4isog}.

Generally, the $z$-component of a degree $\ell$ isogeny onto $y^3=z\,(z-1)$ 
will be a genus 1 Belyi function with the passport $[3^{\ell}/3^{\ell}/3^{\ell}]$,
since the isogeny is an unramified covering. 
The endomorphisms on $y^3=z\,(z-1)$ will be compositions with
genus 0 maps with the passports $[3^k/3^k/3^{k-1}\,1^3]$ 
or $[3^k\,1/3^k\,1/3^k\,1]$ 
with $k=\lfloor \ell/3 \rfloor$.
Recall that a $j=0$ elliptic curve has complex multiplication by $\omega_3$.
Similarly as for $j=1728$ curves, the endomorphisms 
form a ring isomorphic to $\IZ[\omega_3]$. 
The degree of an isogeny equals the squared norm of the corresponding 
algebraic integer in $\IZ[\omega_3]$. As described in \cite[\S 8]{vhpg}, \cite{vhpgell},
these isogenies induce  hypergeometric identities for 
\begin{align*} \label{elin6a}
\hpg{2}{1}{1/2,\,1/6}{7/6}{z}=& \, \frac{z^{-1/6}}{2}
\,\int_{1/\!\sqrt[3]{z}}^{\infty}\;\frac{dX}{\sqrt{X^3-1}}, \\[2pt]
\hpg{2}{1}{1/3,\,2/3}{4/3}{z} = & \, z^{-1/3}\,
\int_{1/\!\sqrt[3]{z}}^\infty\;\frac{dX}{(X^3-1)^{2/3}}
\end{align*}
similar to  (\ref{eq:hpg5}).

\section{Physical Motivations and Prospectus}
\label{sec:conc}

Whilst we have focused on obtaining the dessins in their own right, the main motivation for the particular genus one Belyi maps come from physics.
In the context of gauge theories with toric Calabi-Yau moduli spaces \cite{Feng:2000mi}, it is realized that a bipartite graph on the torus $T^2$ encodes a 4-dimensional supersymmetric gauge theory whose space of vacua is an affine, toric Calabi-Yau threefold.
This is obviously important in string theory (cf.~rapid review in \cite{He:2012js}) and the graph is called a ``brane tiling'',
by drawing the bipartite as a tiling of the doubly-periodic plane.

In particular, there is a classification, arranged by the number 
$N_W$ of black/white nodes in the fundamental region of the tiling, for $N_W$ up to 6 in \cite{Davey:2009bp}.
Here, 
$W$ refers to the superpotential $W$ of the quantum field theory, whose Jacobian determines the relations in the quiver and
the number $2N_W$ corresponds to the number of monomial terms in $W$.
It is therefore a clear open problem of interest to the physics community to obtain the explicit Belyi maps which correspond to these bipartite tilings of $T^2$, considered as dessins on elliptic curves.

\begin{table}[t!]
\hspace{-0.7cm}
\begin{tabular}{|c|l|l|l|l|}\hline
\!Degree\! 
& Belyi Map & Dessin & \!Gauge theory\! 
& Reference \\ \hline \hline
3  
 & $\Phi_0$ in \eqref{eq:ecj0}, \S \ref{sec:mapsj0} 
 & Fig.~\ref{fig:HexDessin}\refpart{b-$1$}; \cite[\!(1.1)]{Davey:2009bp}  
 & $\IC^3$ & \cite[\!(53)]{Jejjala:2010vb}, \cite[\!(2.3)]{Hanany:2011ra}  
\\  
4 
 & $\Psi_0$ in \S \ref{sec:square} 
 & Fig.~\ref{fig:SquareDessin}\refpart{a}; \cite[\!(1.2)]{Davey:2009bp}  
 & conifold $\cC$ & \cite[\!(73)]{Jejjala:2010vb}, \cite[\!(6.1)]{Hanany:2011ra} 
\\ 
6 
 & \eqref{deg6}  
 & Fig.~\ref{fig:HexDessin}\refpart{b-$2$}; \cite[\!(2.1)]{Davey:2009bp} 
 & $\IC^2/\IZ_2\times \IC$ & \cite[\!(89)]{Jejjala:2010vb}, \cite[\!(5.13)]{Hanany:2011ra} 
\\ 
8 
 & $\psi_1$ in \eqref{psi11}
 & Fig.~\ref{fig:SquareDessin}\refpart{b}; \cite[\!(2.5)]{Davey:2009bp}  
 & $\IF_0$(I) &  \cite[\!C.2.3]{Jejjala:2010vb}, \cite[\!(3.8)]{Hanany:2011ra} \\ 
8 
& $\psi_2$ in \eqref{psi12}
 & Fig.~\ref{fig:SquareDessin}\refpart{c}; \cite[\!(2.4)]{Davey:2009bp}   
& $L^{222}$(I) &  \cite[\!(201)]{Jejjala:2010vb}, \cite[\!(6.10)]{Hanany:2011ra}
\\ 
9 
 & $\varphi_0$ in \S \ref{sec:isogj0}
 & Fig.~\ref{fig:HexDessin}\refpart{b-$3$}; \cite[\!(3.2)]{Davey:2009bp}   
 & $\IC^3/\IZ_3$ &  \cite[\!C.1.1]{Jejjala:2010vb}, \cite[\!(5.11)]{Hanany:2011ra}
\\ 
10 
 & \eqref{eq:comp5}
 & App.~\refpart{B}; \cite[\!(3.5)]{Davey:2009bp} 
 & $L^{222}$(II) & \cite[\!(7.2)]{Hanany:2011ra}
\\  
12 & $\varphi_1$ in \eqref{eq:comp1}
 & Fig.~\ref{fig:dessin328}\refpart{c}; \cite[\!(3.28)]{Davey:2009bp} 
 & $dP_3$(I) & \cite[\!(5.3)]{Hanany:2011ra}
\\ 
12 & $\psi_3$ in \eqref{eq:isog1c}
 & Fig.~\ref{fig:SquareDessin}\refpart{d}; \cite[\!(3.26)]{Davey:2009bp} 
 & $L^{333}$(I) &  \cite[\!(202)]{Jejjala:2010vb}
\\ 
12 & conj.~$\psi_3$ in \eqref{eq:isog1c}
 & Fig.~\ref{fig:SquareDessin}\refpart{e}; \cite[\!(3.27)]{Davey:2009bp} 
 & $Y^{3,0}$(I) &
\\ 
12 & $z$-comp. in \eqref{eq:4isog}
 & Fig.~\ref{fig:HexDessin}\refpart{b-$4$}
 & $\IC^2/\IZ_4\times\IC$ &
\\ 
12 & \eqref{eq:comp33} on \eqref{eq:ec1728-2}
 & Fig.~\ref{fig:dessin332}\refpart{c} 
 &  $\IF_0$(II) & \cite[\!(7.1)]{Hanany:2011ra}
\\ 
14 & conj. \eqref{eq:E} & App.~\refpart{E2} & $L^{333}$(II) & \\
14 & conj. \eqref{eq:F} & App.~\refpart{F2} & $dP_3$(III) & \\  
16 & $\psi_1^2$ in \S \ref{sec:square} \refpart{f} 
  & Fig.~\ref{fig:SquareDessin}\refpart{f} 
  & $\cC/(\IZ_2\times\IZ_2)$ & \\ 
16 & in \S \ref{sec:square} \refpart{g} 
  & Fig.~\ref{fig:SquareDessin}\refpart{g} 
  & $Y^{4,0}$ & \\ 
16 & in \S \ref{sec:square} \refpart{h} 
  & Fig.~\ref{fig:SquareDessin}\refpart{h} 
  & $\cC/\IZ_4$ &  \\ 
16 & in \S \ref{sec:square} \refpart{i} 
  & Fig.~\ref{fig:SquareDessin}\refpart{i} 
  & $L^{444}$  & \\ 
16 & conj. \eqref{eq:16-2} & App.~\refpart{I2} & $Z^{3,1}$ & \\
18 & \eqref{eq:deg18a} or \eqref{eq:deg18b}  & App.~\refpart{J} & $dP_3$(IV) & \\ 
20 & \eqref{ec:kl} 
  & Fig.~\ref{fig:SquareDessin}\refpart{k} 
  & $\cC/\IZ_5$(I) &  \\ 
20 & conj.~\eqref{ec:kl}
  & Fig.~\ref{fig:SquareDessin}\refpart{l} 
  & $\cC/\IZ_5$(II)  & \\ 
20 & $\!(X$-comp. in \eqref{eq:isogm2}$)^2\!$ 
  & Fig.~\ref{fig:SquareDessin}\refpart{j} 
  & $\cC/\IZ_5$(III) & \\ 
\hline 
\comment{
12 & \eqref{eq:ec28} & $\varphi_2$ in \eqref{eq:phi28} & Figure \ref{fig:dessin334}\refpart{b}  \\ \hline
12 & \eqref{eq:ec52} & $\varphi_3$ in \eqref{eq:comp33} & Figure \ref{fig:dessin332}\refpart{a}  \\ \hline
12 & \eqref{eq:ec1728-1} & $\varphi_3$ & Figure \ref{fig:dessin332}\refpart{b}  \\ \hline
12 & \eqref{eq:ec1728-2} & $\varphi_3$ & Figure \ref{fig:dessin332}\refpart{c}  \\ \hline
12 & \eqref{eq:isog1c} & $\psi_3$ in \eqref{eq:isog1c} &  Figure \ref{fig:SquareDessin} \refpart{b}, \refpart{c} as Galois orbits \\ \hline
14 & Below \eqref{eq:C} & \eqref{eq:C} & App.~fig.~\refpart{C}\\ \hline
14 & Below \eqref{eq:D} & \eqref{eq:D} & App.~fig.~\refpart{D}\\ \hline
14 & Below \eqref{eq:E} & \eqref{eq:E} & App.~fig.~\refpart{E1} \& \refpart{E2} as a Galois orbit\\ \hline
14 & Below \eqref{eq:F} & \eqref{eq:F} & App.~fig.~\refpart{F1} \& \refpart{F2} as a Galois orbit\\ \hline
14 & Below \eqref{eq:G} & \eqref{eq:G} & App.~ fig.~\refpart{G1} \& \refpart{G2} as a Galois orbit\\ \hline
16 & Below \eqref{eq:16-1} & \eqref{eq:16-1} & App.~fig.~\refpart{H} \\ \hline
16 &  \eqref{ec:16-2} & \eqref{eq:16-2} & App.~fig.~\refpart{I1} \& \refpart{I2} as a Galois orbit \\ \hline
16 & \begin{tabular}{l}
	4-isogenies of  \\
	$E_1$ in \eqref{ec:coni}
	\end{tabular} & 
	 \begin{tabular}{l} 4 maps (f) - (i)\\
	  below \eqref{4iso} 
	  \end{tabular}
	  & Figure \ref{fig:SquareDessin} \refpart{f}--\refpart{i} \\ \hline
18 &  Below \eqref{eq:deg18a} & \eqref{eq:deg18a} & App.~ Fig.~\refpart{J} \\ \hline
20 & Below \eqref{ec:kl} & \eqref{ec:kl} & Figure \ref{fig:SquareDessin} \refpart{k}, \refpart{l} as Galois orbits \\ \hline
24 & 3 curves \eqref{eq:isog34} & \begin{tabular}{l}
	3 maps by \\
	sub.~(\ref{eq:isog28}) into (\ref{eq:phi28}) 
	\end{tabular}
	& Figure \ref{fig:d334isog} \refpart{b}--\refpart{d}
	\\ \hline
}
\end{tabular}
\caption{The list of Belyi maps considered in this paper,
  whose dessins correspond  to consistent gauge theories. 
In the first column, ``conj." means conjugation in the quadratic definition field; 
``comp." means ``component".
\label{t:phys}
}
\end{table}

In the ensuing, we will adhere to the nomenclature of \cite{Davey:2009bp}.
That is, given a dessin, its planar dual graph is a quiver with superpotential, $(Q,W)$.
For any such $(Q,W)$, we have a representation variety $\cM$, which is essentially the algebraic variety corresponding to the path algebra of $Q$ quotient the Jacobian ideal of $W$.
In our context $\cM$ will always be an affine toric Calabi-Yau threefold (CY3).

For example, the dessin which is  the hexagonal bipartite tiling of the doubly period plane corresponds to the trivial CY3 $\IC^3$, as drawn in Figure \ref{fig:HexDessin}\refpart{b-$1$}.
This is none other than the dessin corresponding to the Belyi map
 $\Phi_0=(1+Y)/2$  on the elliptic curve $Y^2 = X^3 +1$ 
 in \eqref{eq:ecj0} and as discussed in detail in \S \ref{sec:mapsj0}.
Other archetypal examples of affine toric CY3 include 
orbifolds of $\IC^3$ by $\IZ_m \times \IZ_n$ with $m,n \in \IZ_{\ge 2}$, 
and complex cones over the del Pezzo surfaces $dP_n$ which are $\IP^2$ blown up at $n=0,1,2,3$ generic points.

%
Our main motivation was to compute the Belyi maps for the dessins in \cite{Davey:2009bp} 
and several other cases corresponding to the gauge theories.
It was realized in \cite{Hanany:2011ra,He:2012xw} that the brane tilings as not merely combinatorial objects which graphically encode the gauge theory data.
Instead, the shape of the $T^2$ specified by the $j$-invariant of the $T^2$ endowed with the complex structure of an elliptic curve plays a crucial role in relating to the so-called R-charges in the quantum field theory.
It is therefore important to obtain the explicit Belyi maps for the dessins.
As shown in our paper, looking for compositions was an easy and efficient way to find a lot of explicit Belyi maps,
and we indeed obtained a wealth Belyi maps leading to interesting  theories.
For example, the Belyi maps for \cite[(3.27), (3.28)]{Davey:2009bp} were all unknown 
in the physics literature and yearn to be computed.
We remark that many theories in the catalogue \cite{Davey:2009bp} 
are marked by ``inc.'', which means that as quantum field theories are ``inconsistent'' 
and need to be dismissed on physical grounds. 
Nevertheless, they establish perfectly well-defined dessins.
Moreover, 
we observe that the property of consistency is not an invariant
of the Galois action of dessins.
This is obvious, for instance, from the dessin examples (E) and (F) in the ensuing table.
The relation between the Galois action and the the physics of the gauge theories remains to be explored.

Table \ref{t:phys} lists all the constructed Belyi maps that correspond to the consistent theories.
Most of the other Belyi maps are relegated  to the Appendix. 
In the table, the gauge theories (quiver representation varieties) which appear in \cite{Davey:2009bp} are marked while several which are unpublished or appear elsewhere are also indicated.
In the table, 
the Belyi maps and dessins are arranged by degree (in the first column).
The second column cross-references explicit expressions for Belyi maps;
the elliptic curves are usually defined immediately around the referenced formula. 
The dessins are referenced both in this paper and in \cite{Davey:2009bp}.
The vacuum moduli space $\cM$ (an affine, toric Calabi-Yau variety of complex dimension three) 
for the graph-dual of the dessins, as bipartite graphs on the doubly-periodic plane, 
is presented in the last column 
signifying the gauge 
theory encoded by the dessins. There are a few Belyi maps, such as the simplest two cases in the first two rows,
which have been constructed in the physics literature and we include them for reference.
The majority of the table, such as $dP_3$(III), $dP_3$(IV), $Y^{4,0}$, etc. are new and have been eagerly awaited.

Regarding the shape of the $T^2$, \cite{Hanany:2011ra} found a counter-example to a conjecture raised in \cite{He:2012xw} that the $j$-invariant of the elliptic curve for the dessin should be the same as if the bipartite graph were isoradially embedded.
While this holds for highly symmetric cases such as the hexagonal or square tilings, for  $L^{222}$(II), it was found to not to hold.
In our table, the $L^{333}$(II) and $dP_3$(III) theories are further counter-examples.
Indeed, locally the dessins are equi-angular in that the edges around each vertex give a regular division because the Belyi covering in local coordinates is a power function by construction: its Taylor series begins with a constant 0, 1 or $\infty$ and the next order is at least a quadratic term.
This is incompatible, except in the most symmetric cases, with
\cite{Jejjala:2010vb} isoradial embeddings.
Furthermore, the $j$-invariant for the dessin also does not seem to be preserve under so-called cluster transformations or Seiberg duality \cite{Hanany:2011bs}, which beckons the question as to the implications of such quiver transformations on the dessins, and indeed, vice versa, the role of Galois conjugation in the quiver.
These are all important further directions.

In summary, using composition of Belyi maps from $\IP^1$, which are univariate rational functions of the projective coordinate on the Riemann sphere, together with covering of the $\IP^1$ by an elliptic curve $E$, we have obtained a wealth of explicit genus 1 dessins.
Furthermore, we can also compose maps from $E$ with isogenies of $E$; this also gives us potentially new genus 1 dessins.


We have of course only begun 
this study as our methods are applicable to any genus 0 dessin, of which there is a small database (cf.~\cite{Bose:2014lea}). In addition to producing genus 1 dessins interesting in their own right, we aim to complete as much as possible the construction of explicit Belyi maps in the catalogue of brane-tilings.
Furthermore, there is an upcoming work \cite{newCat} which significantly expands the known catalogue of \cite{Davey:2009bp} and for which many new Belyi maps and companion elliptic curves need to be computed.
We anticipate a companion paper addressing this in the near future.


\appendix
\section{Appendix: Explicit Dessins}
Here  
we present the dessins of the Belyi map that appear in \S \ref{sec:deg10}, \S \ref{sec:higherd} and some of \S \ref{sec:othercomp}. 
Together with Table \ref{t:phys}, and Figures \ref{fig:dessin332}, \ref{fig:d334isog} and \ref{fig:HexDessin}\refpart{a}, this Appendix constitutes all the genus 1 dessins discussed in the paper.
The dessins presented here are supplemented by the genus 0 dessins with cuts, demonstrating the transition as in Figures \ref{fig:dessin328} and \ref{fig:dessin334}.
This also demonstrates a method of getting genus 1 tilings without necessarily computing the Belyi maps.

The dessins have their associated Belyi maps in the main text in \S\ref{sec:comp};
\comment{
in particular, \refpart{A}, \refpart{B}, \refpart{C},
\refpart{D}, \refpart{E1,E2}, \refpart{F1,F2}, \refpart{G1,G2}, 
\refpart{H}, \refpart{I1,I2},
have their Belyi maps given, respectively,
in \eqref{eq:ellc10a}, \eqref{eq:comp5}, \eqref{eq:C}, \eqref{eq:D},
two Galois orbits of \eqref{eq:E}, two Galois orbits of \eqref{eq:F}, two Galois orbits of \eqref{eq:G}, \eqref{eq:16-1},  two Galois orbitds of 
}
Each dessin is drawn both as a doubly-periodic tiling and as a branched cover over the associated genus 0 dessin and we specify the field of definition.
For reference, we cross-reference the equations where the explicit Belyi map is presented in the text next to each dessin.

As mentioned in \S\ref{sec:deg10}, the genus 1 curve for \refpart{A} is defined over $\IQ$, but the Belyi map is defined over $\IQ(\sqrt{-15})$.
Similarly, the genus 1 curves for 
\refpart{F1}, \refpart{F2} are defined over $\IQ(\sqrt7)$ but the Belyi maps are defined over $\IQ(\sqrt{2\sqrt7-7})$; and the cases \refpart{J}, \refpart{N} have $j=0$ but the dessins are defined over $\IQ(\omega_3)$ or $\IQ(i)$. These dessins are not defined over $\IR$,
but they have certain chess-board asymmetry
induced by patching four complex half-planes into a torus as in Figure
\ref{fig:CubicInv}.
Also, these dessins  \refpart{A}, \refpart{F1}, \refpart{F2}, \refpart{J}, \refpart{N}
are symmetric if the colouring of the vertices is ignored.

\noindent
\begin{picture}(440,158)
\put(0,0){\includegraphics[width=440pt]{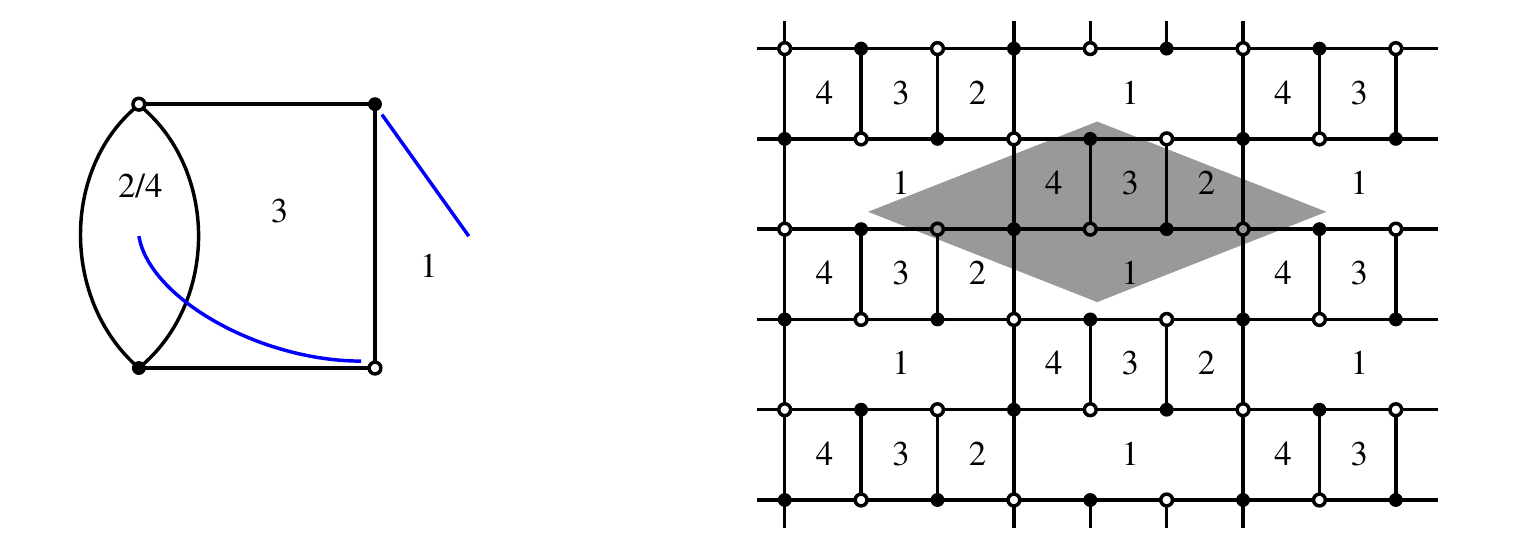}}
\put(0,10){\refpart{A} Map \eqref{eq:deg10a}, $\!j\!=\!-\frac{108}5$;
\cite[\!(3.12)]{Davey:2009bp}.}
\end{picture}
\begin{picture}(440,158)
\put(0,0){\includegraphics[width=440pt]{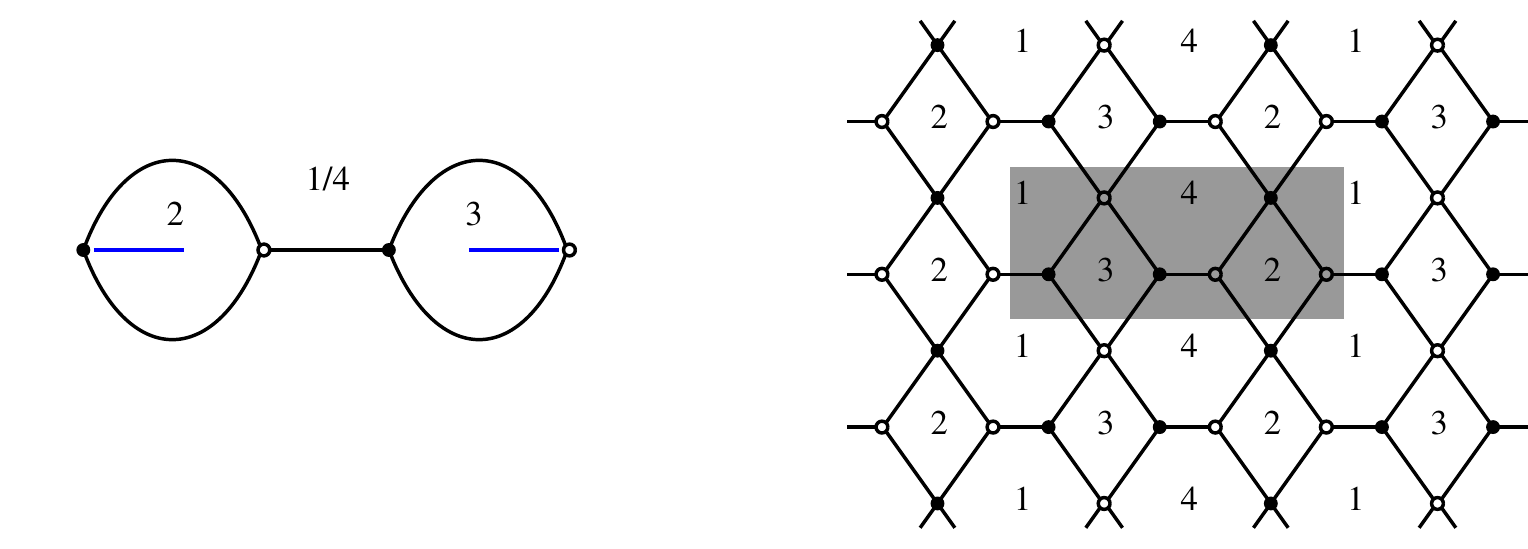}}
\put(0,10){\refpart{B} Map \!\eqref{eq:comp5}, $\!j\!=\!\frac45\,321^3$;\!
\cite[\!(3.5)]{Davey:2009bp}, $\!L^{222}$(II).}
\end{picture}
\begin{picture}(440,158)
\put(0,0){\includegraphics[width=440pt]{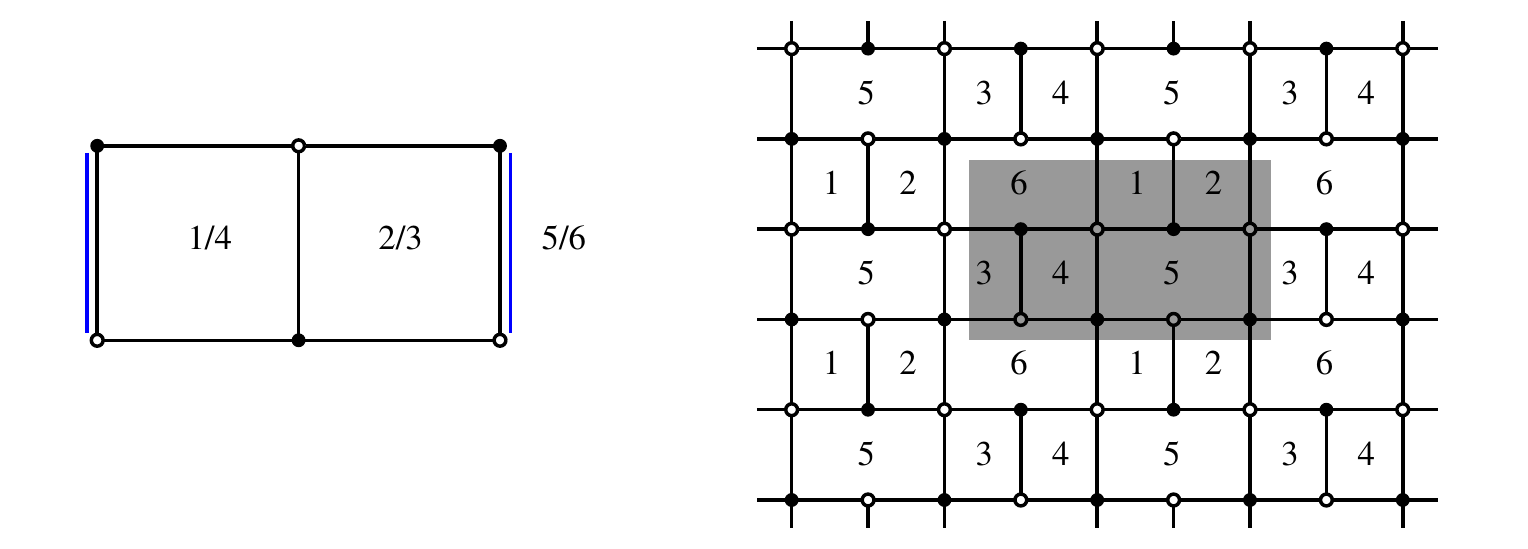}}
\put(0,10){\refpart{C} Map \eqref{eq:C} on a $j=\frac{4}{49}\;57^3$ curve.}
\end{picture}
\begin{picture}(440,158)
\put(0,0){\includegraphics[width=440pt]{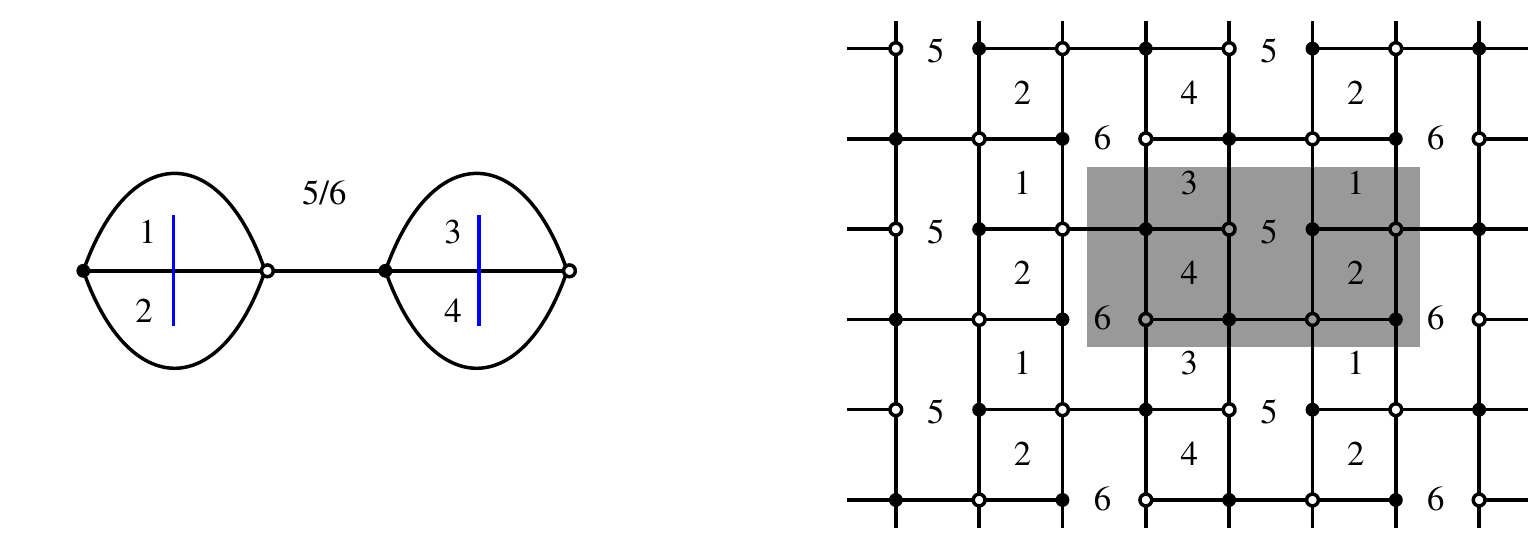}}
\put(0,10){\refpart{D} Map \eqref{eq:D} on a $j=255^3$ curve.}
\end{picture}
\begin{picture}(440,158)
\put(0,0){\includegraphics[width=440pt]{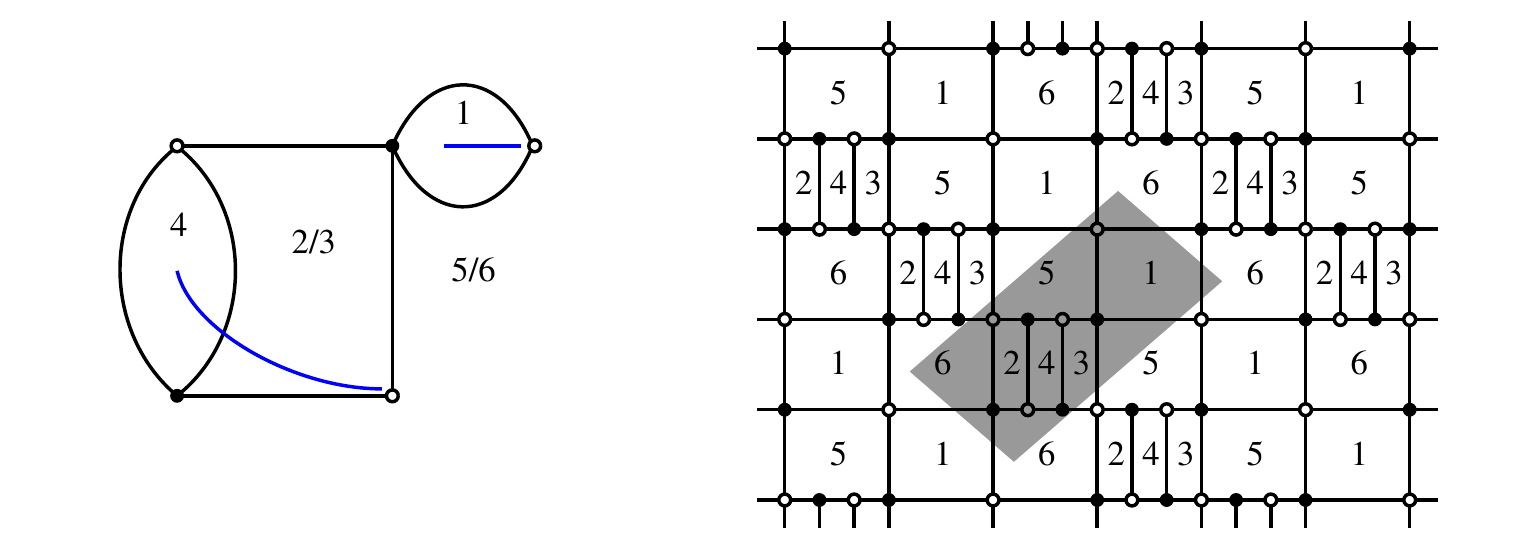}}
\put(0,10){\refpart{E1} Map \eqref{eq:E}; $\xi^3-\xi^2+5\xi+15=0$.}
\end{picture}
\begin{picture}(440,158)
\put(0,0){\includegraphics[width=440pt]{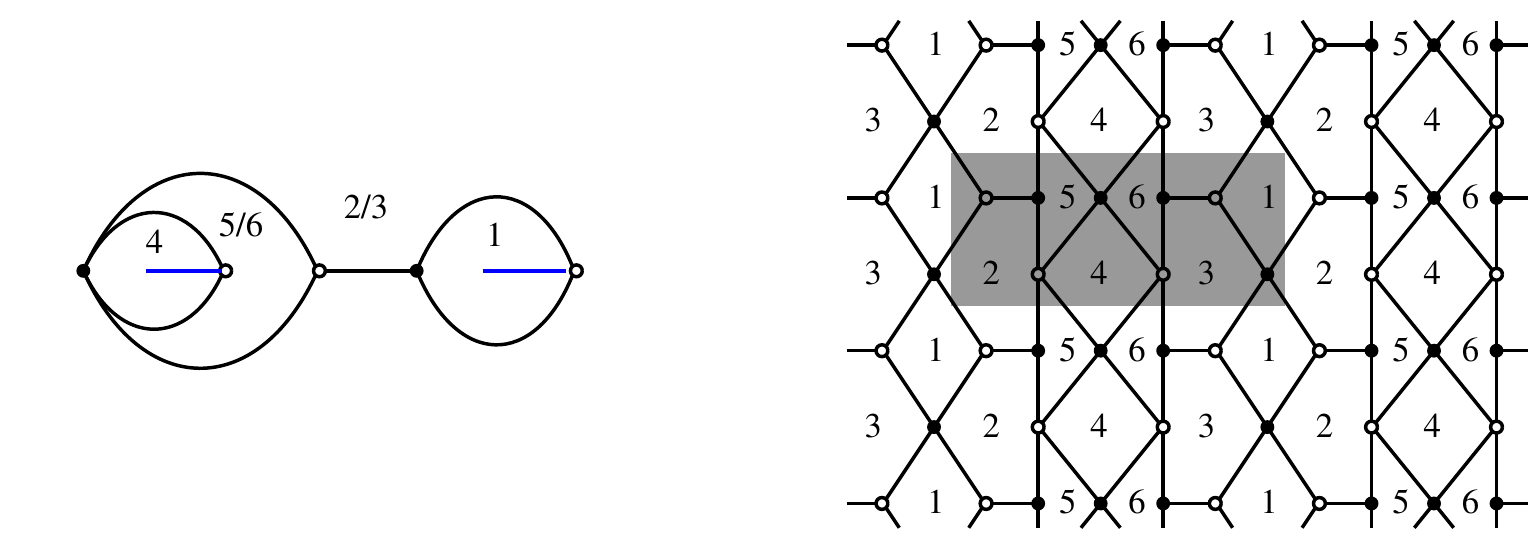}}
\put(0,10){\refpart{E2} Map \eqref{eq:E} with real $\xi$; $L^{333}$(II).}
\end{picture}
\begin{picture}(440,158)
\put(0,0){\includegraphics[width=440pt]{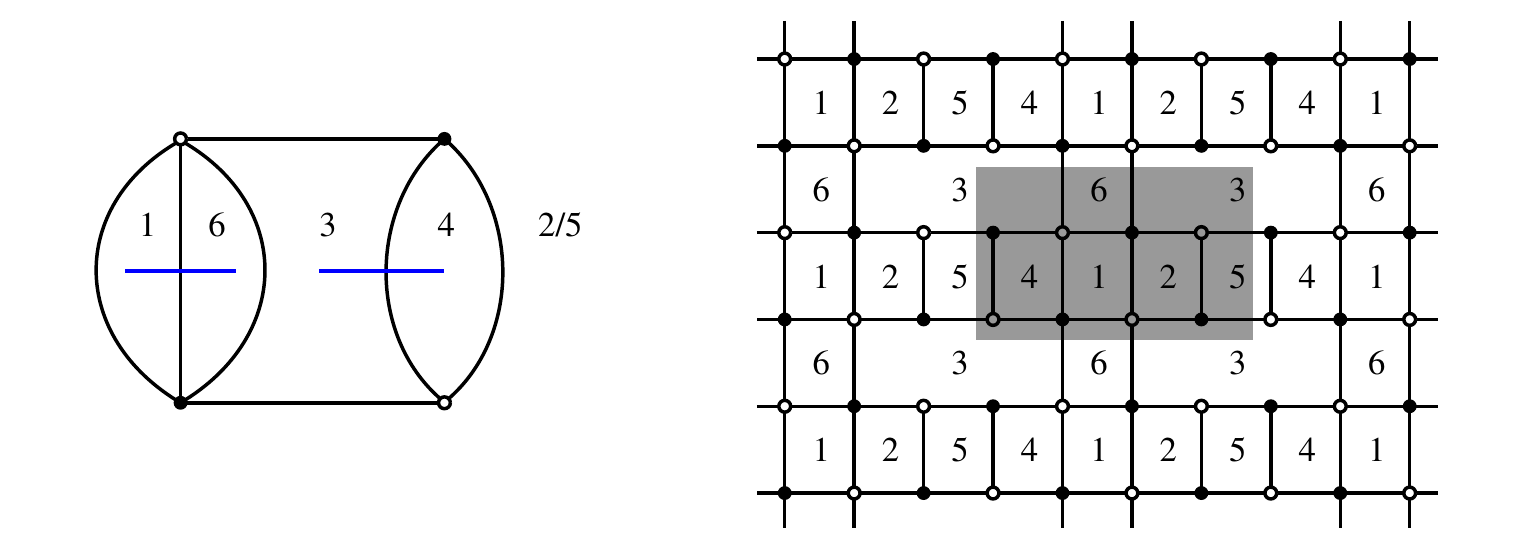}}
\put(0,10){\refpart{F1} Map \eqref{eq:F} on a $\IQ(\sqrt7)$-curve.}
\end{picture}
\begin{picture}(440,158)
\put(0,0){\includegraphics[width=440pt]{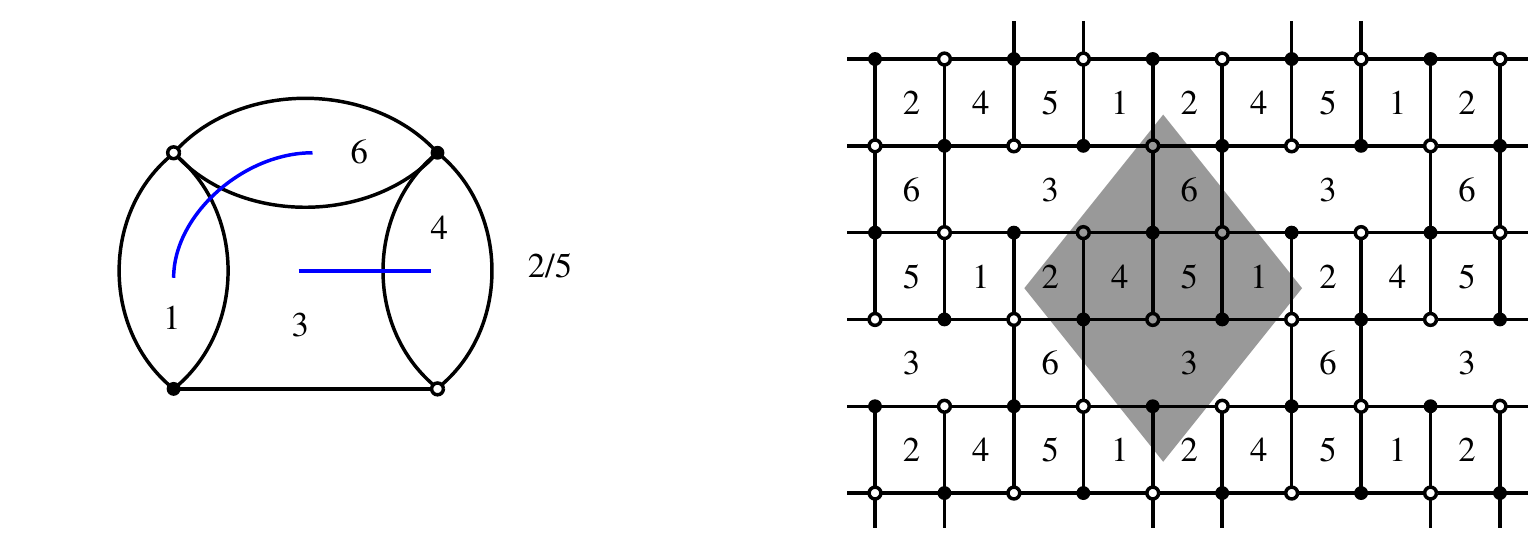}}
\put(0,10){\refpart{F2} $\sqrt7$-conjugated \eqref{eq:F}; $dP_3$(III).}
\end{picture}
\begin{picture}(440,158)
\put(0,0){\includegraphics[width=440pt]{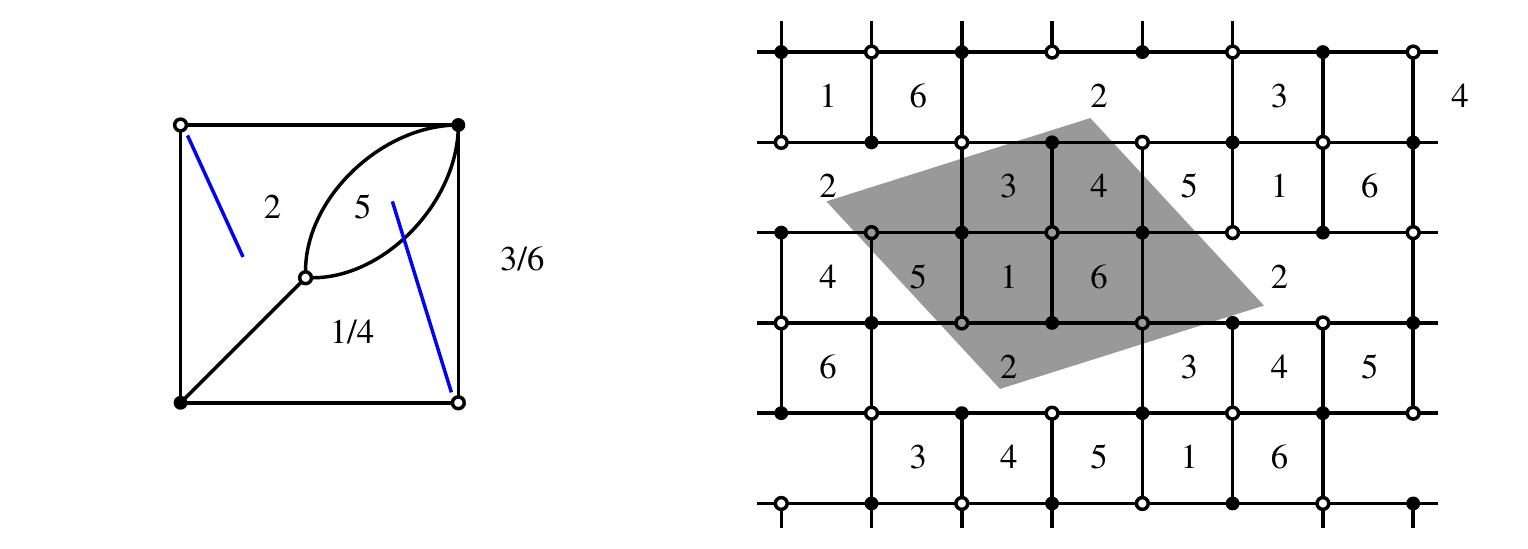}}
\put(0,10){\refpart{G1} Map \eqref{eq:G}; $\eta^3+\eta^2-2\eta+6=0$.}
\end{picture}
\begin{picture}(440,158)
\put(0,0){\includegraphics[width=440pt]{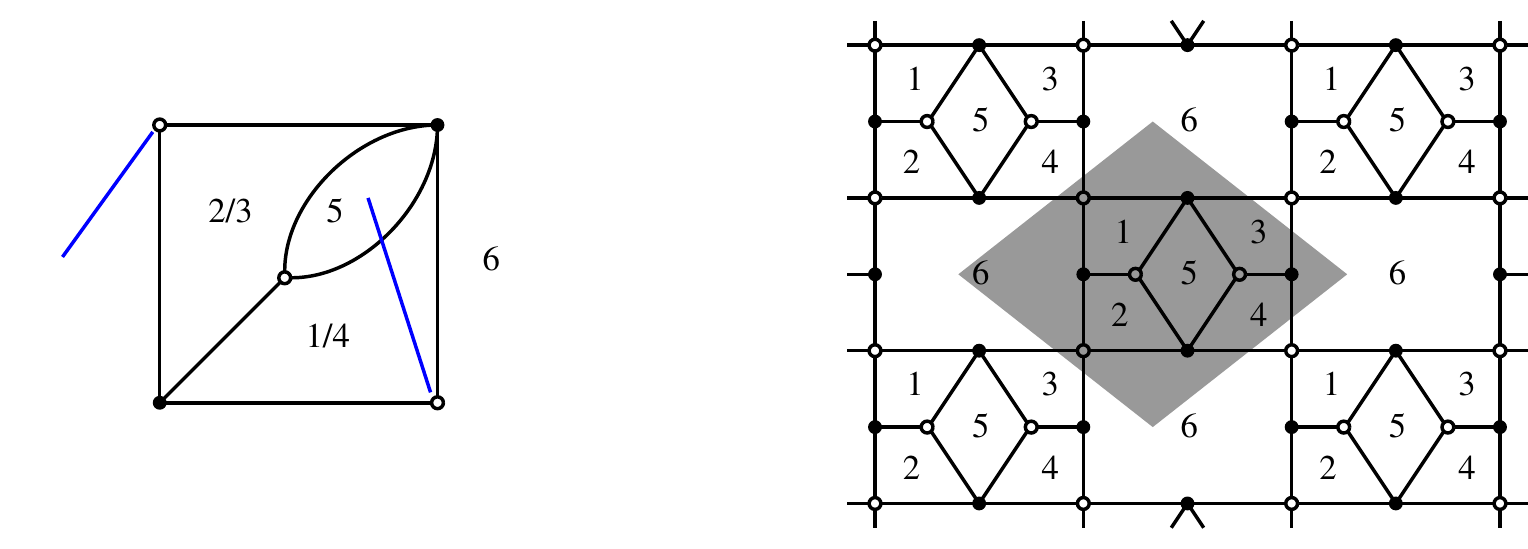}}
\put(0,10){\refpart{G2} Map \eqref{eq:G}, with real $\eta$.}
\end{picture}
\begin{picture}(440,158)
\put(0,0){\includegraphics[width=440pt]{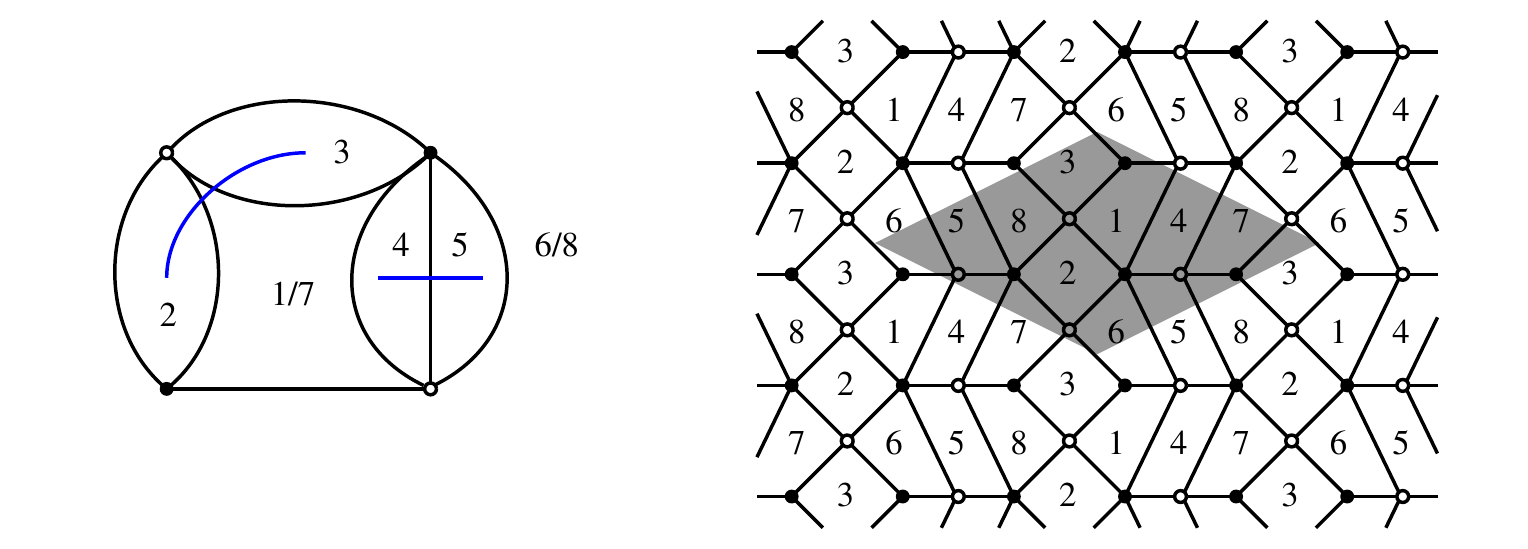}}
\put(0,10){\refpart{H} Map \eqref{eq:16-1} on a $j\!=\!-5000$ curve.}
\end{picture}
\begin{picture}(440,158)
\put(0,0){\includegraphics[width=440pt]{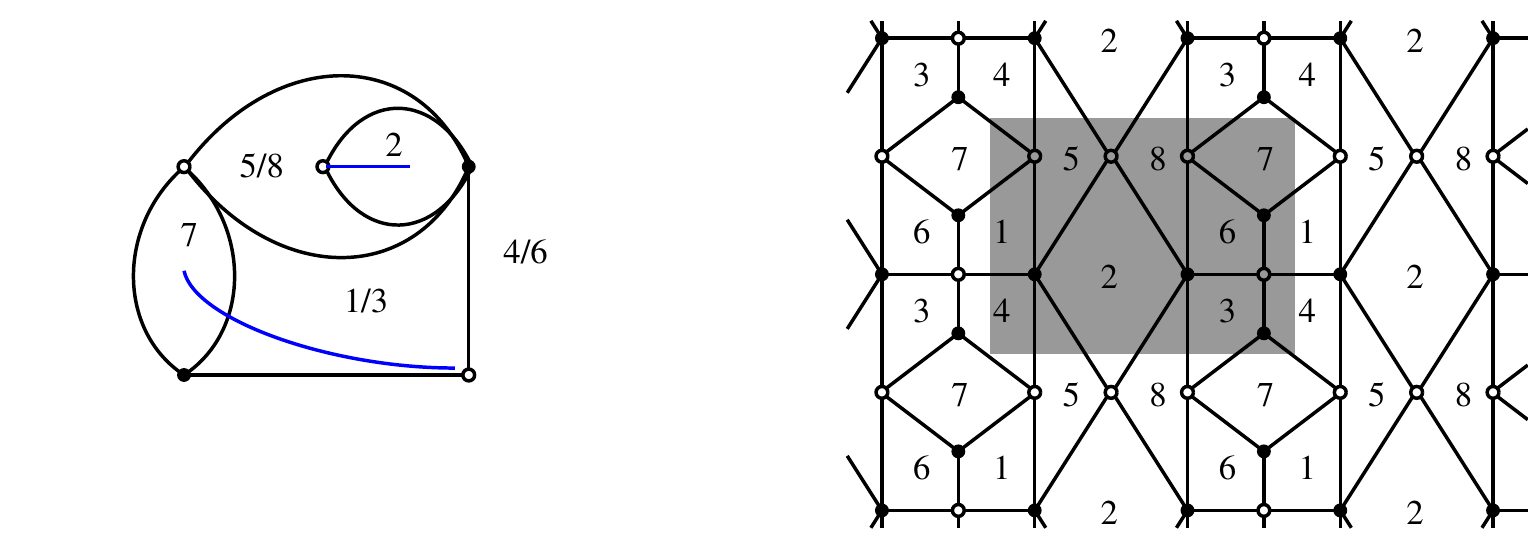}}
\put(0,10){\refpart{I1} Map \eqref{eq:16-2}, defined over $\IQ(\sqrt{10})$.}
\end{picture}
\begin{picture}(440,158)
\put(0,0){\includegraphics[width=440pt]{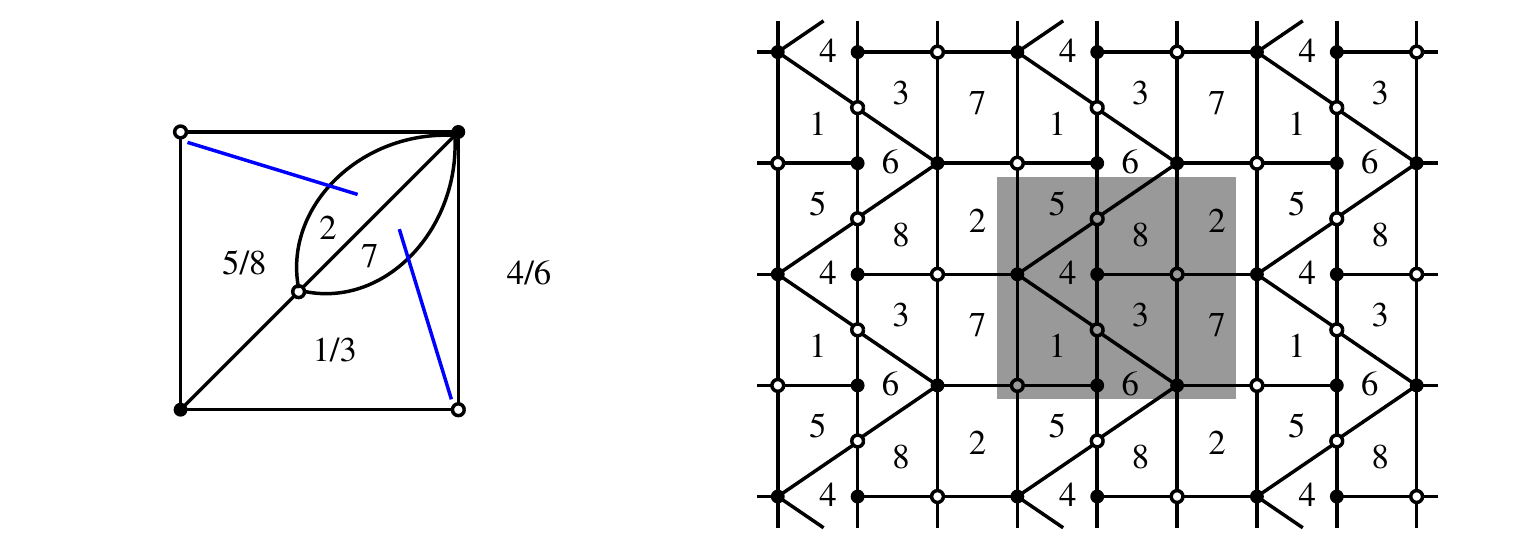}}
\put(0,10){\refpart{I2} $\sqrt{10}$-conjugated \eqref{eq:16-2}; $Z^{3,1}$.}
\end{picture}
\begin{picture}(440,158)
\put(0,0){\includegraphics[width=440pt]{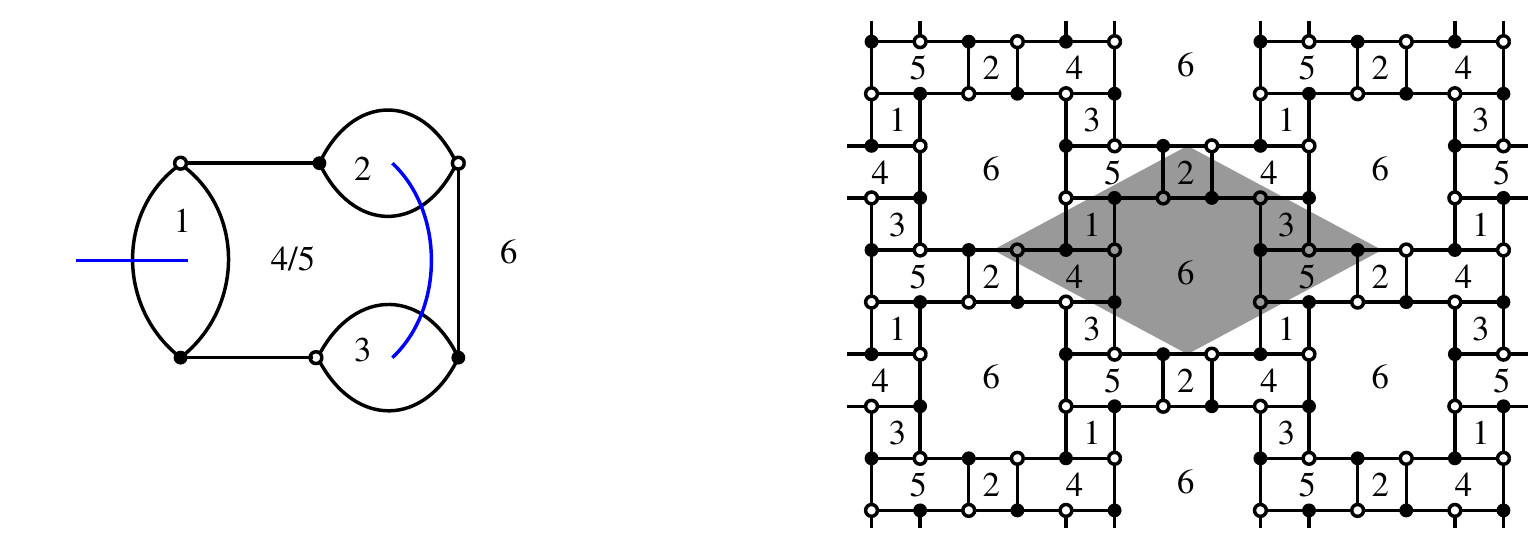}}
\put(0,10){\refpart{J} Map \eqref{eq:deg18a} or \eqref{eq:deg18b}, $j=0$; $dP_3(IV)$.}
\end{picture}
\begin{picture}(440,158)
\put(0,0){\includegraphics[width=440pt]{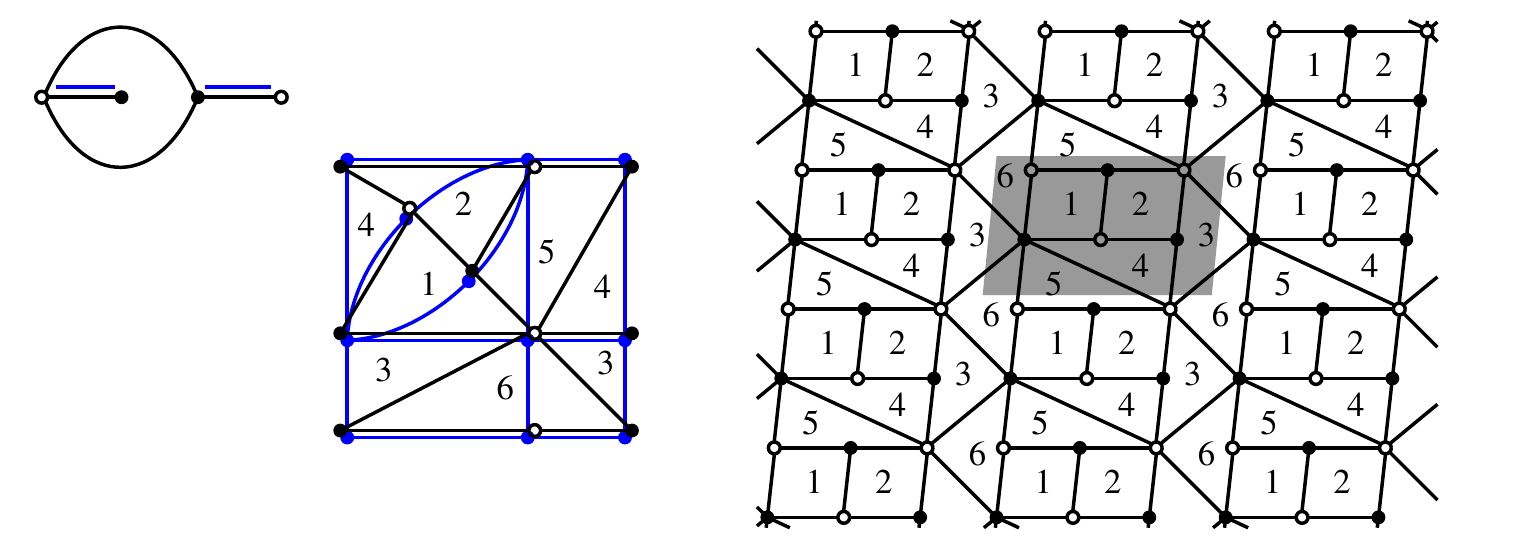}}
\put(0,10){\refpart{K} Map \eqref{eq:phi422} on \eqref{eq:ecK}; \cite[(3.35)]{Davey:2009bp}.}
\end{picture}
\begin{picture}(440,158)
\put(0,0){\includegraphics[width=440pt]{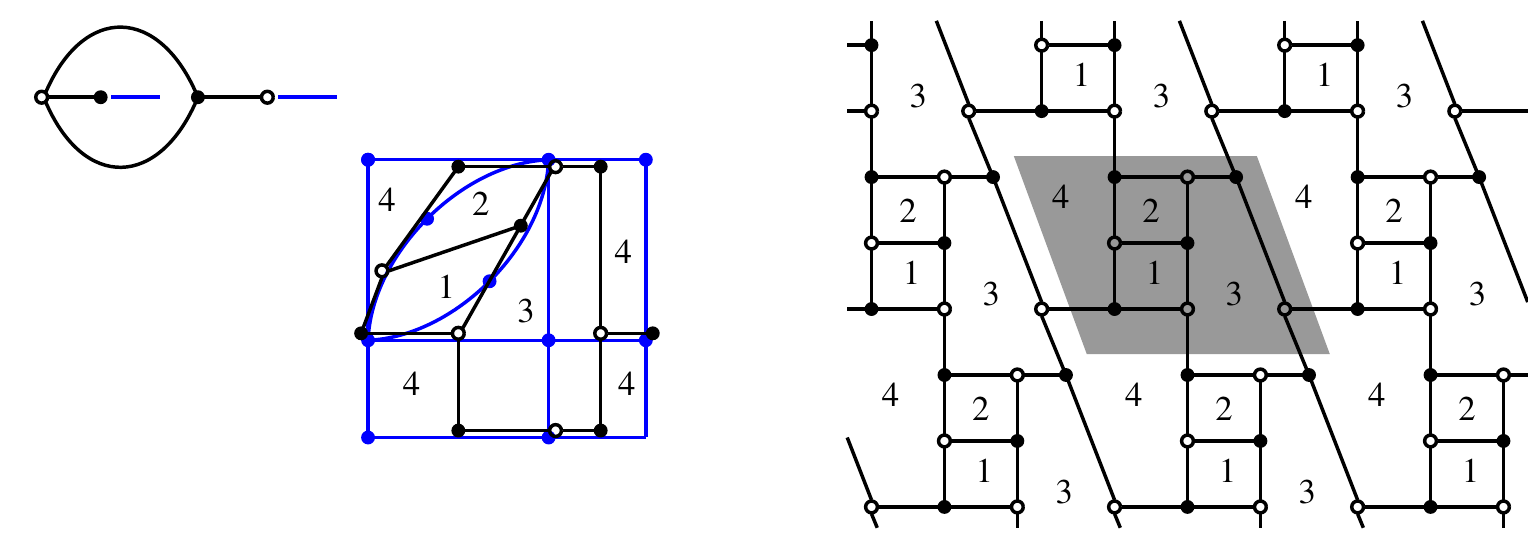}}
\put(0,10){\refpart{L} Map \eqref{eq:phi422} on \eqref{eq:ecL}, with $\sqrt{-2}$.}
\end{picture}
\begin{picture}(440,158)
\put(0,0){\includegraphics[width=440pt]{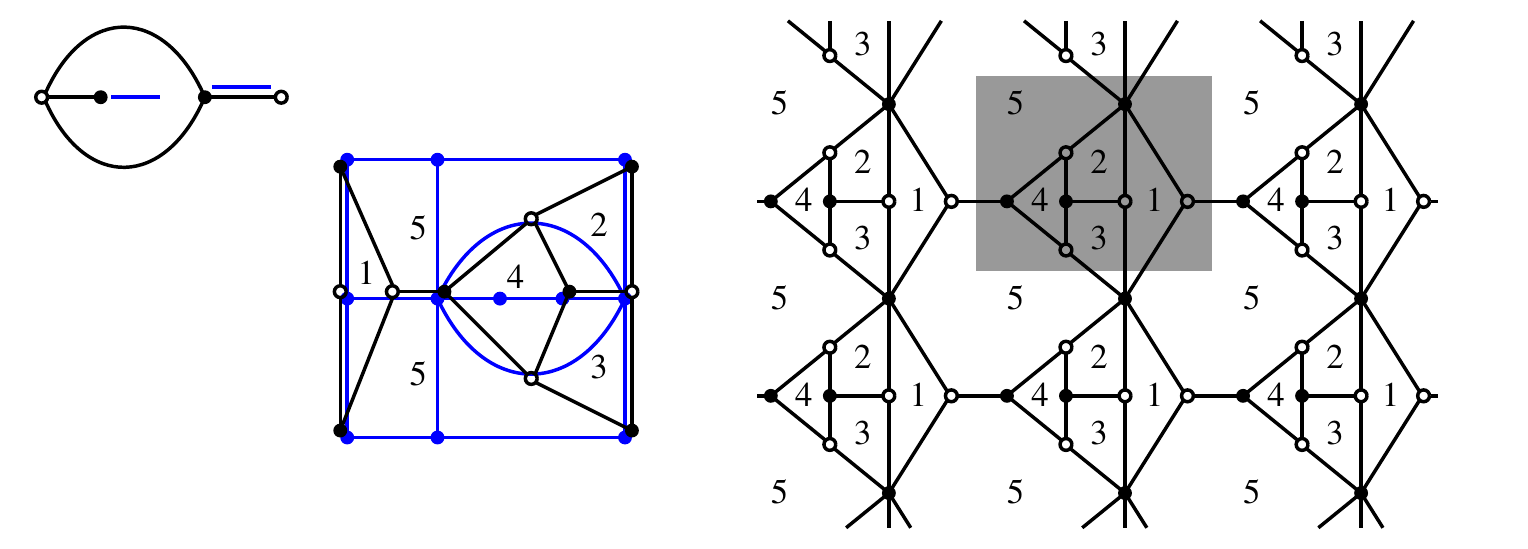}}
\put(0,10){\refpart{M1} Map \eqref{eq:phi422}, with $\sqrt{\sqrt3+1}$.}
\end{picture}
\begin{picture}(440,158)
\put(0,0){\includegraphics[width=440pt]{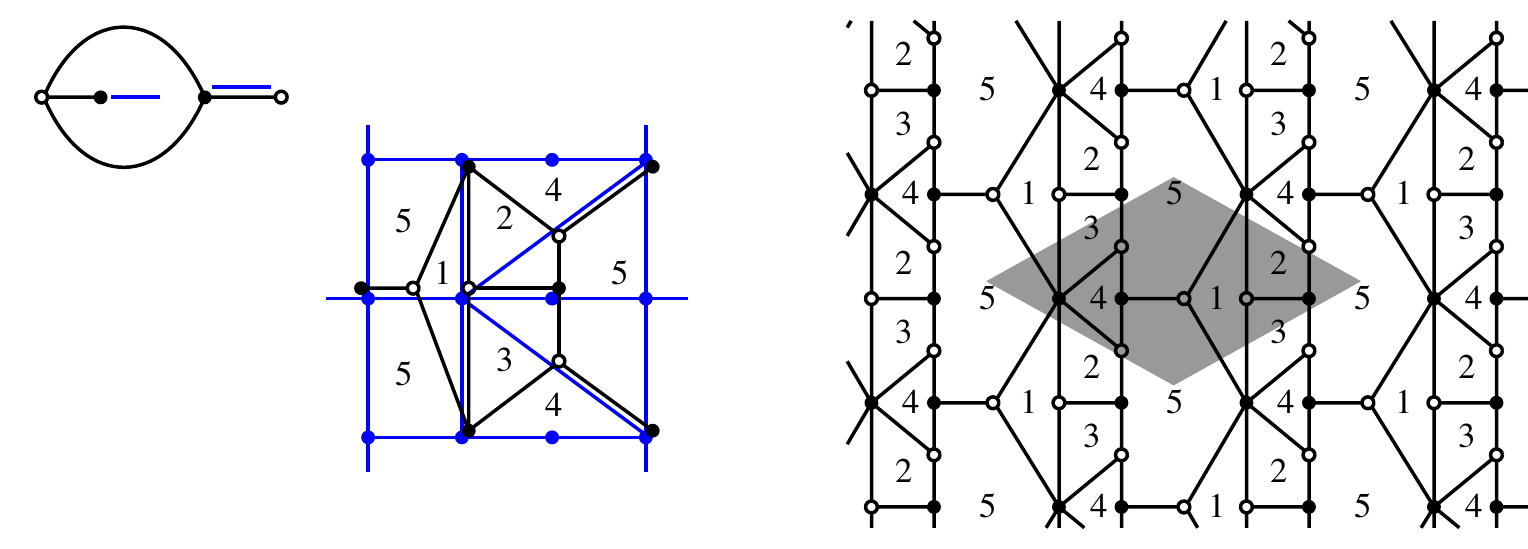}}
\put(0,10){\refpart{M2} Map \eqref{eq:phi422}, with $\sqrt{\sqrt3+1}$.}
\end{picture}
\begin{picture}(440,158)
\put(0,0){\includegraphics[width=440pt]{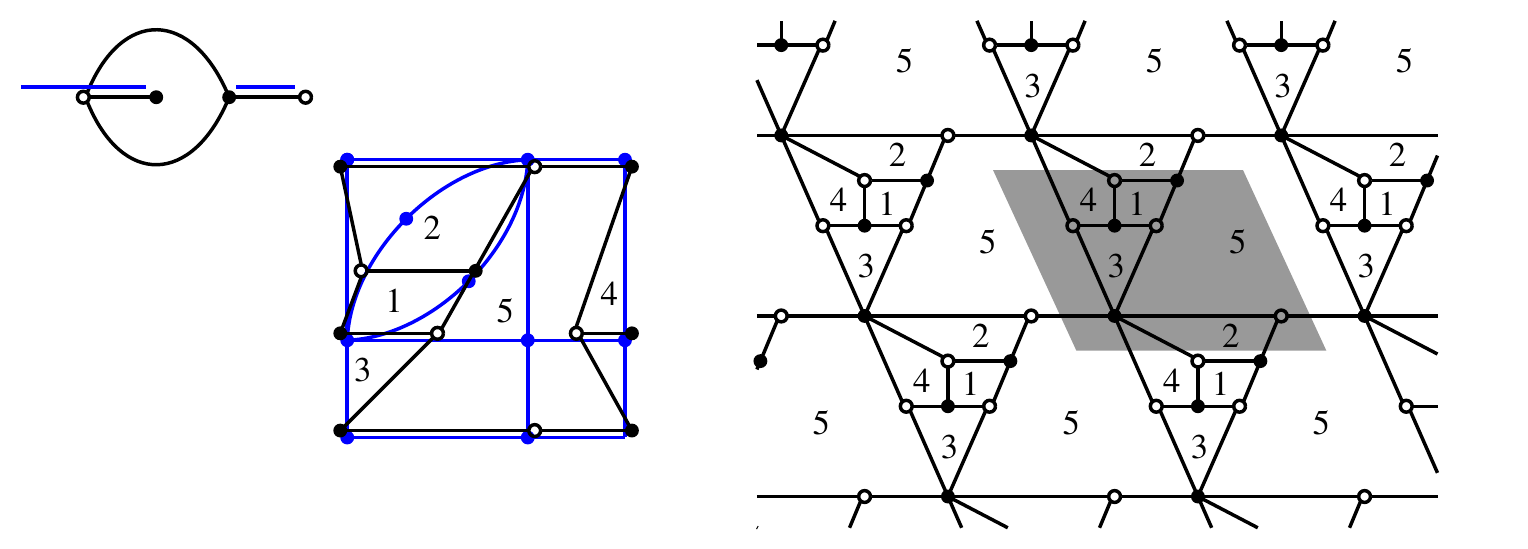}}
\put(0,10){\refpart{M3} Map \eqref{eq:phi422}, with $\sqrt{1-\sqrt3}$.}
\end{picture}
\begin{picture}(440,158)
\put(0,0){\includegraphics[width=440pt]{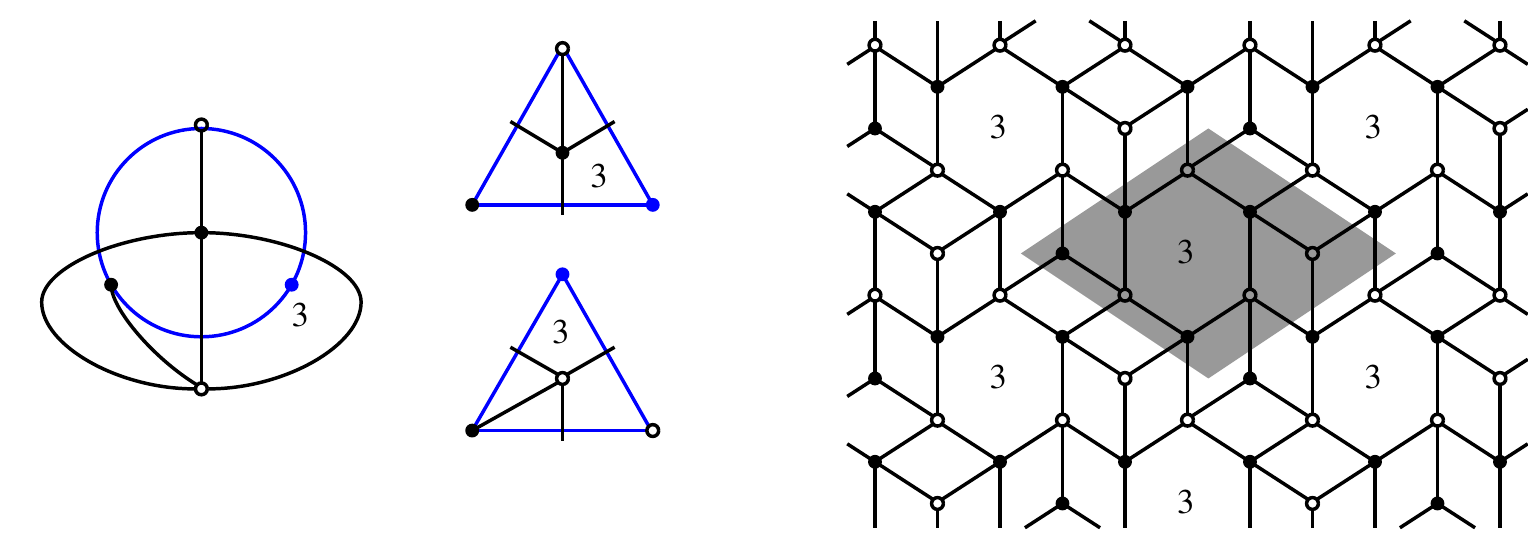}}
\put(0,10){\refpart{N} Map \eqref{eq:deg15} on a $j=0$ curve.}
\end{picture}

\section*{Acknowledgments}
We are most grateful to the organizers for inviting us to the {\it Workshop on Grothendieck-Teichm\"uller Theories}, 2016, at the Chern Institute of Mathematics, Nankai University where the most congenial and productive atmosphere inspired the beginnings of this paper.

YHH would like to thank the Science and Technology Facilities Council, UK, for grant ST/J00037X/1, the Chinese Ministry of Education, for a Chang-Jiang Chair Professorship at NanKai University as well as the City of Tian-Jin for a Qian-Ren Scholarship, and Merton College, Oxford, for her enduring support.


\end{document}